%% file: lineargroupactions.tex
\documentclass[10pt]{amsart}
\usepackage{amsmath}
\usepackage{amssymb}
\usepackage{amsthm}
\usepackage{hyperref}
\input{xcommands}
\title{The Baum-Connes conjecture and proper group actions on affine buildings}
\author{Dmitry Matsnev}
\email{matsnev@math.ist.utl.pt}
\date{\today}
\begin{document}
\address{Departamento de Matem\'atica, Instituto Superior T\'ecnico, Av. Rovisco Pais, 1049-001 Lisboa, Portugal}
\maketitle
\begin{abstract}
We study the possibility of applying a finite-dimensionality argument in order to address parts of the Baum-Connes conjecture for finitely generated linear groups. This gives an alternative approach to the results of Guentner, Higson, and Weinberger concerning the Baum-Connes conjecture for linear groups.\\
For any finitely generated linear group over a field of characteristic zero we construct a proper action on a finite-asymptotic-dimensional $CAT(0)$-space, provided that for such a group its unipotent subgroups have ``bounded composition rank''. The $CAT(0)$-space in our construction is a finite product of symmetric spaces and affine Bruhat-Tits buildings.\\
For the case of finitely generated subgroup of $SL(2,\C)$ the result is sharpened to show that the Baum-Connes assembly map is an isomorphism.
\end{abstract}
\input{introduction}
\input{preliminaries}

\input{asdim}
\input{intchar}
\input{maintheorem}
\input{sl2c}
\section*{Acknowledgements}This paper is based on the author's thesis completed under the supervision of Nigel Higson at Penn State University. The author is very grateful to Professor Higson for his invaluable guidance, comments, and suggestions. 
\bibliographystyle{amsalpha}
\addcontentsline{toc}{chapter}{Bibliography}
\bibliography{bibliography-data}
\end{document}

%% file: xcommands.tex
\usepackage{amsmath,amsthm,amssymb}

\DeclareMathOperator{\asdim}{asdim}
\DeclareMathOperator{\tr}{tr}
\DeclareMathOperator{\Span}{Span}
\DeclareMathOperator{\dist}{dist}

\newcommand{\Q}{\mathbb{Q}}
\newcommand{\A}{\mathbb{A}}
\newcommand{\Z}{\mathbb{Z}}
\newcommand{\C}{\mathbb{C}}
\newcommand{\R}{\mathbb{R}}
\newcommand{\HH}{\mathbb{H}}
\newcommand{\m}[4]{\ensuremath{\begin{pmatrix}{#1}&{#2}\\{#3}&{#4}\end{pmatrix}}}
\newcommand{\h}[3]{\ensuremath{\begin{pmatrix}1&{#1}&{#3}\\0&1&{#2}\\0&0&1\end{pmatrix}}}
\theoremstyle{plain}
\newtheorem{theorem}{Theorem}[section]
\newtheorem*{theorem*}{Theorem}
\newtheorem{lemma}[theorem]{Lemma}
\newtheorem{corollary}[theorem]{Corollary}
\newtheorem*{corollary*}{Corollary}

\theoremstyle{remark}
\newtheorem{remark}[theorem]{Remark}
\theoremstyle{definition}
\newtheorem{definition}[theorem]{Definition}
\newtheorem*{definition*}{Definition}
\theoremstyle{remark}
\newtheorem{example}[theorem]{Example}

%% file: introduction.tex
\section{Introduction}\label{chapterIntroduction}

The Baum-Connes conjecture, introduced in the early 80's by Paul Baum and Alain Connes, connects the $K$-theory of the reduced crossed product of a $C^*$-algebra by a group acting on such algebra and the $K$-homology of the corresponding classifying space of proper actions of that group (for a formal account see~\cite{BaumConnesHigson}).

Let $\Gamma$ be a discrete group acting on a $C^*$-algebra $A$ by automorphisms. The Baum-Connes conjecture proposes that the ``assembly''  morphism 
$$\mu : KK^{\Gamma}(\underline{E}\Gamma,A)\to K(A\rtimes_r\Gamma)$$
from the $K$-homology of the classifying space $\underline{E}\Gamma$ of proper actions of $\Gamma$ to the $K$-theory of the reduced crossed product of $A$ by $\Gamma$ is in fact an isomorphism.

While the Conjecture is formulated in terms of pairs $(\Gamma,A)$, it is possible to state it purely in terms of group $\Gamma$. Namely one can ask whether the original conjecture holds for the group $\Gamma$ and any $C^*$-algebra $A$, on which $\Gamma$ acts. This version of the conjecture is called \emph{the Baum-Connes conjecture with coefficients}\footnote{Some authors refer to it as to the Baum-Connes \emph{property} with coefficients, in the view of the counterexamples (modulo a statement due to Gromov) by Higson, Lafforgue, and Skandalis in~\cite{HigsonLafforgueSkandalis}.}, while the version with $A=\C$ is called \emph{the Baum-Connes conjecture without coefficients}. In this work our main concern will be the conjecture with coefficients. The conjecture with coefficients is stable under passing to subgroups (see~\cite{ChabertEchterhoff}).

Sometimes the actual conjecture is being split into two ingredients: the injectivity and surjectivity of the assembly map $\mu$. We shall be mostly concerned with injectivity.
In fact, even this part of the conjecture is a rather difficult problem. It is known, for example, that the injectivity of the Baum-Connes assembly map implies the Novikov's higher signature conjecture~\cite{FerryRanickiRosenberg}.

There are different approaches to the understanding of the Baum-Connes conjecture. Here we will touch only few of them, directly referring to the notions of proper action and asymptotic dimension.

\begin{theorem}[Guentner, Higson, and Weinberger,~\cite{GuentnerHigsonWeinberger}]\label{theoremGL2Std}
For any field $K$, the Baum-Connes conjecture with coefficients holds for any countable subgroup of $GL(2,K)$.
\end{theorem}

\begin{theorem}[Guentner, Higson, and Weinberger,~\cite{GuentnerHigsonWeinberger}]\label{theoremGLnStd}
For any field $K$ and any natural number $n$, the injectivity portion of the Baum-Connes conjecture with coefficients holds for any countable subgroup of $GL(n,K)$.
\end{theorem}

The original proofs of these theorems were based on the existence of  metrically proper isometric actions of such groups on a Hilbert space (and appealing to the results of Higson and Kasparov in~\cite{HigsonKasparov}) and uniform embeddings of these groups  into a Hilbert space (and thus referring to results of Yu in~\cite{Yu2000}).

Unfortunately, in these proofs the group actions were not constructed explicitly, and Hilbert spaces on which the groups act can easily be infinite-dimensional.

There is another approach, which reflects the coarse geometric view of discrete groups. It is based on the following result:

\begin{theorem}[Yu,~\cite{Yu}]\label{theoremYu}
If a finitely generated group $\Gamma$ has finite asymptotic dimension as a discrete metric space with the word-length metric corresponding to a finite set of generators and its classifying space is of finite homotopy type, then the coarse Baum-Connes conjecture holds for $\Gamma$, or, alternatively, the Baum-Connes assembly map for such $\Gamma$ is injective.
\end{theorem}

This theorem was proven in~\cite{Yu} by means of extensive $\epsilon$-$\delta$ calculations and is based on the elementary facts from the asymptotic dimension theory only. Higson in~\cite{Higson} proved that the finiteness assumption on the classifying space can be relaxed. There exists an alternative proof by Wright in~\cite{Wright}, which is based on the introduction of a more delicate coarse structure on a metric space.

We are interested in investigating the existence of proper actions of linear groups on finite-dimensional spaces in the spirit of  Theorem~\ref{theoremYu}. Also we would like to give a more elementary proof of Theorem~\ref{theoremGL2Std}, without appealing to Hilbert space techniques.

As a result, we prove the following theorem:

\begin{theorem}\label{theoremMain}
Let $K$ be any field of characteristic 0, and $\Gamma$ be a finitely generated subgroup of $SL(n,K)$. If there exists a natural number $m$, such that all unipotent subgroups of $\Gamma$ have composition rank bounded by $m$ (informally, this condition means that they are built out of ``layers'' of asymptotic dimension not more than $m$), 
then $\Gamma$ acts properly on a finite-asymptotic-dimensional $CAT(0)$ space.
\end{theorem}

For the case $n=2$ the construction is sharpened:

\begin{theorem}\label{theoremMain2}
Let $K$ be any field of characteristic 0, and $\Gamma$ be a finitely generated subgroup of $SL(2,K)$. Then $\Gamma$ satisfies the Baum-Connes conjecture (with coefficients).
\end{theorem}



Note that since $GL(n,K)$ embeds into $SL(n+1,K)$, any finitely generated subgroup of  $GL(n,K)$ can be realized as a finitely generated subgroup of $SL(n+1,K)$, for which we know how to construct a proper action. Thus our result extends to all finitely generated linear groups over $K$.

For the case of $2\times2$ matrices, we also can extend our result to any finitely generated subgroup $\Gamma$ of $GL(2,K)$. Consider the following exact sequence:
\begin{equation*}
1\to\Gamma\cap SL(2,K)\to\Gamma\to K^{\times},
\end{equation*}
where the only nontrivial maps are the natural inclusion and the determinant respectively. Our argument proves the conjecture for $\Gamma\cap SL(2,K)$, the group $K^{\times}$ is abelian and so is any subgroup of it, and we can apply the extension result of Chabert, Echterhoff, and OyonoOyono from~\cite{ChabertEchterhoffOyonoOyono} to conclude that $\Gamma$ itself satisfies the Baum-Connes conjecture with coefficients.

As a subproduct of our discussion, we also compute asymptotic dimension of affine buildings corresponding to linear groups over discretely valued fields; asymptotic dimension of symmetric spaces was computed by Bell and Dranishnikov in~\cite{BellDranishnikov04}.  

%% file: preliminaries.tex
\section{Affine buildings associated to discrete valuations}\label{sectionPreliminaries}

In this section we collect some technical tools to be used later in the discussion.

\subsection{Discrete valuations}\label{sectionPreliminariesDiscreteValuations}

\begin{definition}\label{definitionDiscreteValuation}
Let $R$ be an integral domain (a commutative ring without zero divisors in which $0\ne1$). A map $\nu:R\to\Z\cup\{+\infty\}$ is called a \emph{discrete valuation} if it satisfies the following properties for any $a, b\in R$:
\begin{itemize}
\item $\nu(a)=+\infty$ if and only if $a=0$
\item $\nu(ab)=\nu(a)+\nu(b)$
\item $\nu(a+b)\geq\min(\nu(a),\nu(b))$
\end{itemize}
\end{definition}

Given a discrete valuation $\nu$ of an integral domain $R$, it can be extended to the field of fractions  frac$(R)$ by
\begin{equation*}
\nu\left(\frac ab\right)=\nu(a)-\nu(b),\qquad\qquad a, b\in R.
\end{equation*}

Starting from a field $K$, one can regard $K$ as a ring and define the corresponding notion of a discrete valuation, following Definition~\ref{definitionDiscreteValuation}.

To any discrete valuation $\nu$ one can associate its \emph{ring of integers}, usually denoted by $\mathcal O_{\nu}$. It consists of all elements of the field $K$ with non-negative valuation.

Any element $\pi$ of $\mathcal O_{\nu}$ with valuation $1$ is called a \emph{uniformizer} of the  valuation $\nu$.

Note that the ideal $\pi\mathcal O_{\nu}$ is maximal in $\mathcal O_{\nu}$~\cite{Cassels}. The \emph{residue field} of the valuation $\nu$ is the quotient of $\mathcal O_{\nu}$ by $\pi\mathcal O_{\nu}$.

\begin{example}
Fix a prime number $p$ and define a $p$-adic valuation $\nu_p$ on $\mathbb Q$ by
$$\nu_p\left(p^n\frac ab\right)=n, a, b, n\in\mathbb Z \mbox{ and } (a,p)=(b,p)=1.$$
The ring of integers of $\nu_p$ consists of all rational numbers without any occurrence of $p$ in the denominator, $p$ itself serves as a uniformizer, and the residue field is isomorphic to $\Z/p\Z$.
\end{example}

\begin{example}
Another valuation which we shall employ is a discrete valuation on $K[x]$ given by $\nu(f)= -\deg f$. The extension of this valuation to the field of fractions $K(x)$ is given by
\begin{equation*}
\nu\left(\frac{q(x)}{r(x)}\right)= \deg r-\deg q, \mbox{ where }q, r\in K[x].
\end{equation*}
The ring of integers here consists of all ratios of polynomials $\frac qr$ with $\deg q\leqslant \deg r$, the residue field in this case is infinite.
\end{example}

\subsection{Affine buildings for $SL(n,k)$}\label{sectionSLbuildings}\label{sectionPreliminariesAffineBuildings}

The notion of a \emph{building} as a special kind of a simplicial complex was developed by Bruhat and Tits in the early 1970's (see~\cite{Tits}). We shall use a certain type of buildings, namely affine buildings, associated to the special linear group over a field with discrete valuation, the construction and basic properties of which are briefly outlined here. Refer to~\cite{Garrett} for more details.

Suppose $K$ is a field with a discrete valuation $\nu$ and denote by $\mathcal O$ the ring of integers of $\nu$ and by $\pi$ some uniformizer.

Let $\mathcal L$ be the set of $\mathcal O$-lattices in $K^n$, and let $\bar{\mathcal L}$ be the set of homothety classes of lattices, i.e. $\bar{\mathcal L}=\mathcal L/K^*$. One constructs an $(n-1)$-dimensional simplicial complex, the Bruhat-Tits building $X$, in the following way:

\begin{itemize}
\item The vertices of $X$ are the elements of $\bar{\mathcal L}$.
\item A set of $(m+1)$ distinct vertices $\bar L_0, \bar L_1,\dots,\bar L_m$ forms an $m$-simplex in $X$ if
\begin{equation*}
L_0\subset L_1\subset\dots\subset L_m\subset\pi^{-1}L_0
\end{equation*}
for some representatives $L_0, L_1,\dots,L_m$ of $\bar L_0, \bar L_1,\dots,\bar L_m$ respectively.
\end{itemize}

To shorten notation we shall write $L_{m+1}$ for $\pi^{-1}L_0$ in the chain of inclusions above.

In particular, any maximal simplex in $X$, a \emph{chamber}, is given by a chain of lattices
\begin{equation}\label{FirstChamber}
L_0\subset L_1\subset\dots\subset L_n\subset\pi^{-1}L_0=L_{n+1},
\end{equation}
such that for any $i$ between $0$ and $m$ the quotient $L_{i+1}/L_i$ is a one-dimensional vector space over $\varkappa=\mathcal O/\pi\mathcal O$, the residue field of the valuation $\nu$.
 
This chamber has $n$ subsimplices of dimension $(n-2)$, \emph{facets}, each obtained by omitting one lattice in the chain above, say the $i^{th}$ facet is given by
$$L_0\subset L_1\subset\dots\subset L_{i-1}\subset L_{i+1}\subset\dots\subset L_n\subset\pi^{-1}L_0=L_{n+1}.$$
Any other chamber sharing this facet with the former one is given by
\begin{equation}\label{SecondChamber}
L_0\subset L_1\subset\dots\subset L_{i-1}\subset L\subset L_{i+1}\subset\dots\subset L_n\subset L_{n+1},
\end{equation}
where $L/L_{i-1}$ is a one-dimensional subspace of the two-dimensional $\varkappa$-vector space $L_{i+1}/L_{i-1}$.

\begin{example}
For the case $n=2$ the affine building described above is an infinite simplicial tree with chambers being the edges of the tree and the facets of these edges being the vertices of the tree. The valence of each vertex is equal to the cardinlity of the residue field (and can be infinite).
\end{example}

In general, one may regard an affine building as a generalization of a simplicial tree as a  $CAT(0)$-space of nonpositive curvature for higher dimensions.

The natural action of $SL(n,K)$ on $K^n$ gives rise to an action of $SL(n,K)$ on $\bar{\mathcal L}$ and therefore to a simplicial action on $X$. This enables us to talk about the action of any subgroup of $SL(n,K)$ on $X$.

We finish this discussion with the following simple observation.

\begin{remark}
The stabilizer of a vertex $x\in X$ is a subgroup of $SL(n,K)$ conjugate in $SL(n,K)$ to $SL(n,\mathcal O)$. This means, in particular, that the coefficients of the characteristic polynomial of any element of such stabilizer are in $\mathcal O$.
\end{remark}

%% file: asdim.tex
\section{Asymptotic dimension}\label{chapterAsdim}

The notion of \emph{asymptotic dimension} of a metric space was introduced by Gromov  in~\cite{Gromov} as an asymptotic property of a metric space. We shall use the following definition taken from~\cite{RoeLecturesOnCoarseGeometry}.

\begin{definition}\label{definitionAsdim}
A metric space $X$ is said to have asymptotic dimension not more than $n$ if for any $d>0$ there exists a uniformly bounded cover of $X$ with the property that any ball of radius $d$ in $X$ meets not more than $(n+1)$ elements of the cover.
\end{definition}

The minimal bound $n$ from the definition above is called the \emph{asymptotic dimension} of $X$ and is denoted by $\asdim X$. 

Now we record a few simple properties of this notion.

\begin{itemize}
\item If $X\subseteq Y$, then $\asdim X\leqslant\asdim Y$.
\item $\asdim(X\cup Y)=\max\{\asdim X,\asdim Y\}$.
\item $\asdim(X\times Y)\leqslant\asdim X+\asdim Y$.
\end{itemize}

Given a finitely generated group $\Gamma$, one may regard it as a (discrete) metric space, namely, fix a symmetrized finite set of generators $S=S^{-1}$ and define the word-length metric on $\Gamma$ to be
\begin{equation}\label{definitionWordMetric}
\dist(\gamma_1,\gamma_2)=\mbox{length of a shortest word in } S \mbox{ representing }\gamma_1^{-1}\gamma_2.
\end{equation}

Different metrics for the same group, which arise this way, are Lipschitz equivalent. Bearing this in mind, one defines the asymptotic dimension of a finitely generated group to be the asymptotic dimension of the underlying discrete metric space.

Note that the metric, defined via~\eqref{definitionWordMetric}, restricts to a metric on any subset of $\Gamma$, in particular, one can equip any subgroup of $\Gamma$ (not necessarily finitely generated) with such a compatible metric.

\begin{definition}\label{definitionAsdimGroup}
A (not necessarily finitely generated) group $G$ has asymptotic dimension not more than $n$ if for any finitely generated subgroup $H$ of $G$, $\asdim H\leqslant n$. The infimum of all such $n$ over all possible finitely generated subgroups of $G$ is called the asymptotic dimension of $G$.
\end{definition}

Definition~\ref{definitionAsdimGroup} is consistent in the following sense. If $H$ is a finitely generated subgroup of a finitely generated group $G$, then $\asdim H\leqslant\asdim G$. Indeed, the word metrics, induced on $H$ by its own set of generators and by the restriction of the corresponding metric on $G$, coming from its set of generators, are coarsely equivalent. This means that the ``intrinsic'' $\asdim H$ and the one, calculated using the metric on $G$, coincide. 

To complete this discussion, we recall the notion of a proper group action.

\begin{definition}
An action of a metrizable group $\Gamma$ on a metrizable space $X$ is called \emph{proper} if there exist finite subgroups $H_i$ of $G$, $i\in I$, and the $H_i$-invariant subspaces $X_i$ of $X$ such that the natural maps $p_i:G\times_{H_i}X_i\to X$ are $G$-homeomorphisms onto their images with
\begin{equation*}
\bigcup\limits_{i\in I}p_i(G\times_{H_i}X_i)=X.
\end{equation*}
\end{definition}

Along with a notion of a proper action, there exists a notion of a metrically proper action.

\begin{definition}
An isometric action of a discrete group $G$ on a metric space $X$ is called \emph{metrically proper} if for all bounded subsets $Y\subseteq X$ the set
\begin{equation*}
\left\lbrace \left. g\in G \right| g.Y\cap Y\ne\emptyset \right\rbrace 
\end{equation*}
is finite.
\end{definition}

\begin{remark}
If a discrete group $G$ admits a metrically proper action on a space $X$ of finite asymptotic dimension, then $G$ has finite asymptotic dimension in its own metric, and its dimension does not exceed the one of $X$. To see this note that in the setup above $G$ is coarsely equivalent to the orbit $Gx_0$, $x_0\in X$, because of properness, while the latter is a subspace of $X$, and so its asymptotic dimension does not exceed the one of $X$.
\end{remark}

\subsection{Example: affine buildings}\label{sectionAsdimUltrametric}

Recall that a metric space $X$ is called \emph{ultrametric} if its distance function $\dist(\cdot,\cdot)$ satisfies the following strong triangle property:
\begin{equation}\label{equationUltrametricTriangle}
\dist(x,z)\leq\max(\dist(x,y),\dist(y,z)) \qquad\qquad\mbox{ for any } x, y, z\in X.
\end{equation}

The last condition in Definition~\ref{definitionDiscreteValuation} of a discrete valuation corresponds to condition~\eqref{equationUltrametricTriangle},
which means that a discretely valued field is an ultrametric space with regard to the associated metric. Now we shall calculate its asymptotic dimension.

\begin{lemma}\label{theoremAsdimUltrametric}
For an ultrametric space $X$ its asymptotic dimension is $0$.
\end{lemma}
\begin{proof}
Suppose $d>0$ is chosen. We shall construct a cover of $X$ with the property that any $d$-ball lies entirely in one of the elements of the cover.

According to the ultrametric condition~\eqref{equationUltrametricTriangle} the closed $d$-ball centered at $y$ has the property that any two points $x, z$ inside such a ball are not more than distance $d$ apart. This means that any point not in this ball has distance strictly more than $d$ from any point inside the ball, and, moreover, any two $d$-balls with different centers are at least $d$-disjoint.

Thus the collection of all $d$-balls in $X$ is a $d$-disjoint cover of $X$ with required properties.
\end{proof}

Our estimate of the asymptotic dimension of affine buildings corresponding to discrete valuations will be based on the following result of Bell and Dranishnikov:

\begin{theorem}[{cf.~\cite[Theorem~1]{BellDranishnikov04}}]\label{theoremAsdimLipschitz}
Let $f:X\to Y$ be a Lipschitz map of a (geodesic) metric space to a metric space. If for every $R>0$ the inverse images of the $R$-balls $f^{-1}(B_R(y))$ are isometric for all $y\in Y$ and have asymptotic dimension $m$, then $\asdim X\leqslant\asdim Y+m$.
\end{theorem}

Throughout this section we fix a discretely valued field $K$ with valuation $\nu$ and denote by $X$ the building associated to $SL(n,K)$.

We shall show that $X$ is coarsely equivalent to a certain matrix group equipped with a pseudometric coming from a length function and thus reduce the computation of the asymptotic dimension of $X$ to the one of that matrix group. The basic definition follows.

\begin{definition}\label{definitionLengthFunction}
A \emph{length function} on a group $G$ is function $l:G\to[0,\infty)$ such that
\begin{enumerate}
\item $l(1)=0$,
\item $l(g)=l(g^{-1})$, and
\item $l(gh)\leqslant l(g)+l(h)$
\end{enumerate}
for any $g, h$ in $G$.
\end{definition}

Following~\cite{GuentnerHigsonWeinberger}, introduce a length function on $SL(n,K)$ in the following way:
\begin{equation}\label{equationLengthOnSLn}
l(g)=-\min_{1\leqslant i,j\leqslant n}\left\lbrace \nu(g_{ij}), \nu(g^{ij})\right\rbrace, \qquad g\in SL(n,K).
\end{equation}
Here we use the standard notation: $g_{ij}$ is the $(i,j)$-th entry of matrix $g$, and $g^{ij}$ is the $(i,j)$-th entry of its inverse. This length function defines a pseudometric on $SL(n,K)$ and on its subgroups by
\begin{equation}\label{equationPseudometric}
dist(g,h)=l(g^{-1}h), \qquad g,h\in SL(n,K).
\end{equation}

We start with some technical computations. Let $G$ be the group $SL(n,K)$ and let $B$ denote its upper-triangular subgroup:
\begin{equation*}
B=\left\lbrace g\in G \; | \; g_{ij}=0 \mbox{ for } 1\leqslant j<i\leqslant n \right\rbrace.
\end{equation*}
Further, let $N$ be the uni-upper-triangular subgroup of $B$, that is its elements satisfy $g_{ii}=1$ for all $i=1,\dots,n$. Let $A$ denote the subgroup of diagonal matrices in $B$. All these subgroups inherit the pseudometric~\eqref{equationPseudometric} from $G$.

\begin{lemma}\label{lemmaInequalities}
Let $a$ and $b$ be two real numbers. Then
\begin{equation*}
-\min\left\lbrace 0, a, \frac12b \right\rbrace \leqslant -\min\left\lbrace 0, a, b \right\rbrace \leqslant -2\min\left\lbrace 0, a, \frac12b \right\rbrace.
\end{equation*}
\end{lemma}

\begin{lemma}\label{lemmaAsdimN}
$\asdim N=0$.
\end{lemma}
\begin{proof}
Define a length function $\tilde l$ on $N$ in the following fashion. For $g\in N$ let
\begin{equation}\label{equationTildeLength}
\tilde l(g)=-\min_{1\leqslant i<j\leqslant n}\left\lbrace 0, \frac1{2^{j-i-1}}\nu(g_{ij}), \frac1{2^{j-i-1}}\nu(g^{ij})) \right\rbrace.
\end{equation}

Before going further, we mention that $\tilde l$ is just a slight modification of the length function $l$ in~\eqref{equationLengthOnSLn}, where the effect of each matrix entry diminishes exponentially as we move away from the diagonal. In fact, for $2\times2$ matrices $\tilde l$ and $l$ coincide. To see this, recall that $\nu(1)=0$ and $\nu(0)=+\infty$. For any element $g=\m{1}{z}{0}{1}$ of $N$ we have
\begin{equation*}
l(g)=-\min\left\lbrace \nu(1), \nu(z), \nu(0), \nu(-z) \right\rbrace=-\min\left\lbrace 0, \nu(z) \right\rbrace=\tilde l(g).
\end{equation*}
Note that we do not claim that $\tilde l$ defines a length function on the entire $G$, but rather on its subgroup $N$ only.

We shall show that $\tilde l$ is coarsely equivalent to the restriction of the original length function $l$ onto $N$ and that $N$ becomes an ultrametric space with respect to $\tilde l$. This will allow us to appeal to Theorem~\ref{theoremAsdimUltrametric}.

To simplify notation, we shall concentrate on the case of $3\times3$ matrices; the general case is similar. 

We start by checking that $\tilde l$ is indeed a length function on $N$. The only nontrivial condition from Definition~\ref{definitionLengthFunction} which $\tilde l$ does not satisfy right away is the last one. Take two elements of $N$:
\begin{equation*}
g=\h{x}{y}{z}\quad\mbox{ and }\quad h=\h{u}{v}{w}, \; x,y,z,u,v,w\in K.
\end{equation*}
For convenience, let us rewrite~\eqref{equationTildeLength} for the case of $g$, $h$, and their product explicitly:
\begin{gather*}
\tilde l(g)=-\min\left\lbrace 0, \nu(x), \nu(y), \frac12\nu(z), \frac12\nu(xy-z)\right\rbrace, \\
\tilde l(h)=-\min\left\lbrace 0, \nu(u), \nu(v), \frac12\nu(w), \frac12\nu(uv-w)\right\rbrace, \mbox{ and}\\
\tilde l(gh)=-\min\left\lbrace 0,\nu(x+u),\nu(y+v),\frac12\nu(xv+z+w),\frac12\nu(xy+uv+yu-z-w)\right\rbrace.
\end{gather*}

We shall show that \emph{each} of the quantities
\begin{equation}\label{equationQuantitiesList}
0,\nu(x+u),\nu(y+v),\frac12\nu(xv+z+w), \mbox{ and } \frac12\nu(xy+uv+yu-z-w)
\end{equation}
is greater or equal to
\begin{equation}\label{equationMinList}
m=\min\left\lbrace 0, \nu(x), \nu(y), \frac12\nu(z), \frac12\nu(xy-z),\nu(u), \nu(v), \frac12\nu(w), \frac12\nu(uv-w)\right\rbrace.
\end{equation}
It is sufficient to show that \emph{each} quantity in~\eqref{equationQuantitiesList} is greater or equal to \emph{some} quantity from the list in~\eqref{equationMinList}.
\begin{enumerate}
\item $0\geqslant0$.
\item $\nu(x+u)\geqslant\min\left\lbrace\nu(x),\nu(u)\right\rbrace$ in accordance with the additive property of a discrete valuation. Obviously $\min\left\lbrace\nu(x),\nu(u)\right\rbrace\geqslant m$.
\item $\nu(y+v)\geqslant m$ in a similar fashion.
\item $\frac12\nu(xv+z+w)\geqslant\frac12\min\left\lbrace\nu(xv),\nu(z),\nu(w)\right\rbrace$ to start with. Note that both $\frac12\nu(z)$ and $\frac12\nu(w)$ occur in ~\eqref{equationMinList}, and we only need to check that $\frac12\nu(xv)\geqslant m$. If the last inequality were not true, then $\frac12\nu(xv)$ would be strictly less than any element from the list in~\eqref{equationMinList}, in particular, $\nu(x)$ and $\nu(v)$. But this gives us a contradiction, for the average $\frac12\nu(xv)=\frac12(\nu(x)+\nu(v))$ can not be strictly less 
than both of its components $\nu(x)$ and $\nu(v)$.
\item $\frac12\nu(xy+uv+yu-z-w)\geqslant\frac12\min\left\lbrace\nu(xy),\nu(uv),\nu(yu),\nu(z),\nu(w)\right\rbrace$ in accordance with the additive property of a discrete valuation.  Now $\frac12\nu(xy)=\frac12(\nu(x)+\nu(y))$, and this quantity is greater or equal at least one of $\nu(x)$ and $\nu(y)$ from the list in~\eqref{equationMinList}. Similarly $\frac12\nu(uv)$ and $\frac12\nu(yu)$ are greater or equal at least one quantity among $\nu(u), \nu(v)$, and $\nu(y)$ respectively. Valuations of $z$ and $w$ occur in~\eqref{equationMinList} themselves, so that we even do not need to do any arithmetics for them.
\end{enumerate}

The statement that we have just proven can be written as
\begin{gather*}
\min\left\lbrace 0,\nu(x+u),\nu(y+v),\frac12\nu(xv+z+w),\frac12\nu(xy+uv+yu-z-w)\right\rbrace \geqslant\\
\min\lbrace \min\left\lbrace 0, \nu(x), \nu(y), \frac12\nu(z), \frac12\nu(xy-z)\right\rbrace ,\\
\min\left\lbrace 0, \nu(u), \nu(v), \frac12\nu(w), \frac12\nu(uv-w)\right\rbrace\rbrace.
\end{gather*}
Clearly this implies
\begin{equation*}
\tilde l(gh)\leqslant\max\left\lbrace \tilde l(g), \tilde l(h)\right\rbrace.
\end{equation*}

Finally, we want to show that $l$ and $\tilde l$ are coarsely equivalent on $N$. In fact, for any $g\in N$ we have $\tilde l(g)\leqslant l(g)\leqslant2\tilde l(g)$. To see this, let
\begin{equation*}
g=\h{x}{y}{z}, \quad x,y,z\in K,
\end{equation*}
so that
\begin{gather*}
l(g)=-\min\left\lbrace 0, \nu(x), \nu(y), \nu(z), \nu(xy-z)\right\rbrace \quad\mbox { and}\\
\tilde l(g)=-\min\left\lbrace 0, \nu(x), \nu(y), \frac12\nu(z), \frac12\nu(xy-z)\right\rbrace.
\end{gather*}
An invocation of Lemma~\ref{lemmaInequalities} with $a=\min\left\lbrace\nu(x),\nu(y)\right\rbrace$ and $b=\min\left\lbrace\nu(z),\nu(xy-z)\right\rbrace$ gives us the desired inequalities. Note that for the general case of $n\times n$ matrices one needs to use Lemma~\ref{lemmaInequalities} repeatedly $(n-2)$ times.
\end{proof}

\begin{lemma}\label{lemmaAsdimA}
$\asdim A=n-1$.
\end{lemma}
\begin{proof}
Take an element $g$ in $A$ and define its ``coarse version'' $\tilde g$ by letting $\tilde g_{ij}=\pi^{\nu(g_{ij})}$ (here $\pi$ denotes a uniformizer of the valuation $\nu$). We have:
\begin{equation*}
dist(g,\tilde g)=l(g^{-1}\tilde g)=-\min_{1\leqslant i\leqslant n}\left\lbrace \nu((g_{ii})^{-1}\tilde g_{ii}), \nu((\tilde g_{ii})^{-1}g_{ii}) \right\rbrace=0.
\end{equation*}
Thus $A$ is coarsely equivalent to the set of all such $\tilde g$, that is, to its own subgroup
\begin{equation*}
\tilde A=\left\lbrace\left. \begin{pmatrix}\pi^{\nu_1} & \cdots & 0 & 0 \\ \vdots & \ddots & \vdots & \vdots \\ 0 & \cdots & \pi^{\nu_{n-1}} & 0 \\ 0 & \cdots & 0 & \pi^{-\nu_1-\dots-\nu_{n-1}}\end{pmatrix} \; \right| \; \nu_1,\dots,\nu_{n-1}\in\Z \right\rbrace.
\end{equation*}
This subgroup, however, is isomorphic to $\Z^{n-1}$ (the matrix shown above corresponds to a point $(\nu_1,\dots,\nu_{n-1})$ in $\Z^{n-1}$) with the length function
\begin{equation*}
\tilde l((\nu_1,\dots,\nu_{n-1}))=\max\left\lbrace |\nu_1|,\dots,|\nu_{n-1}|,|\nu_1+\dots+\nu_{n-1}|\right\rbrace.
\end{equation*}
This length function is bounded from below by the $l^{\infty}$ length function and from above by the $l^1$ length function; for a finite-dimensional space all these functions are Lipschitz equivalent and furnish a metric space of dimension $(n-1)$ (see~\cite{Gromov} for details).

It is known that asymptotic dimension is a coarse invariant (see~\cite{Gromov}), therefore $\asdim A=n-1$ as well.
\end{proof}

Now we are in a position to state

\begin{theorem}
The affine building $X$ has finite asymptotic dimension.
\end{theorem}
\begin{proof}
First we note that our building $X$ is coarsely equivalent to the group $G=SL(n,K)$. Then if $C$ is the maximal compact subgroup of $G$, one can represent $G$ as $CB$, where $B$ consists of upper-triangular matrices, so that $G$ itself is coarsely equivalent to $B$.

Since asymptotic dimension is a coarse invariant, it will be sufficient to show that $B$ has finite asymptotic dimension in order to prove that $X$ does.

Define a homomorphism $f:B\to A$ in the following fashion:
\begin{equation*}
f: \begin{pmatrix}g_{11} & g_{12} & \cdots & g_{1n} \\ 0 & g_{22} & \cdots & g_{2n} \\ \vdots & \vdots & \ddots & \vdots \\ 0 & 0 & \cdots & g_{nn}\end{pmatrix}\mapsto\begin{pmatrix}g_{11} & 0 & \cdots & 0 \\ 0 & g_{22} & \cdots & 0 \\ \vdots & \vdots & \ddots & \vdots \\ 0 & 0 & \cdots & g_{nn}\end{pmatrix}.
\end{equation*}
We claim that $f$ is a $1$-Lipschitz map. Since $f$ is a homomorphism, we need to compare $l(f(g))$ and $l(g)$ for $g\in B$. According to~\eqref{equationLengthOnSLn},
\begin{equation*}
l(g)=-\min\limits_{1\leqslant i\leqslant j\leqslant n}\left\lbrace \nu(g_{ij}), \nu(g^{ij}) \right\rbrace,
\end{equation*}
whilst
\begin{equation*}
l(f(g))=-\min\limits_{1\leqslant i=j\leqslant n}\left\lbrace \nu(g_{ij}), \nu(g^{ij}) \right\rbrace,
\end{equation*}
which is obviously not more than the former quantity.

Given an $R>0$, the $R$-ball in $A$ centered at $a$ is
\begin{equation*}
B_R(a)=\left\lbrace\left. g\in A \right| |\nu(g_{ii})-\nu(a_{ii})|\leqslant R, i=1,\dots,n \right\rbrace.
\end{equation*}
Therefore
\begin{equation*}
f^{-1}\left( B_R(a)\right) =\left\lbrace\left. g\in B \right| |\nu(g_{ii})-\nu(a_{ii})|\leqslant R, i=1,\dots,n \right\rbrace.
\end{equation*}
Left multiplication by $a^{-1}$ is an isometry between $f^{-1}\left( B_R(a)\right)$ and the inverse of the $R$-ball centered at $1\in A$. The latter set is
\begin{gather*}
\left\lbrace\left. g\in B \right| |\nu(g_{ii})|\leqslant R, i=1,\dots,n \right\rbrace=
\bigcup_{\stackrel{g_{ii}\in K}{|\nu(g_{ii})|\leqslant R}}
\begin{pmatrix}g_{11} & \cdots & 0 \\ \vdots & \ddots & \vdots \\ 0 & \cdots & g_{nn}\end{pmatrix}N=\\
\bigcup_{\stackrel{\nu_1,\dots,\nu_n\in\Z}{-R\leqslant \nu_1,\dots,\nu_n\leqslant R}}\begin{pmatrix}\pi^{\nu_1} & \cdots & 0 \\ \vdots & \ddots & \vdots \\ 0 & \cdots & \pi^{\nu_n}\end{pmatrix}\bigcup_{\stackrel{g_{ii}\in K}{\nu(g_{ii})=0}}\begin{pmatrix}g_{11} & \cdots & 0 \\ \vdots & \ddots & \vdots \\ 0 & \cdots & g_{nn}\end{pmatrix}N.
\end{gather*}
We claim that the subgroup
\begin{equation*}
\tilde N=\bigcup_{\stackrel{g_{ii}\in K}{\nu(g_{ii})=0}}\begin{pmatrix}g_{11} & \cdots & 0 \\ \vdots & \ddots & \vdots \\ 0 & \cdots & g_{nn}\end{pmatrix}N
\end{equation*}
is coarsely equivalent to $N$. Indeed, for any element of $\tilde N$, say 
\begin{equation*}
g'=\begin{pmatrix}g_{11} & \cdots & 0 \\ \vdots & \ddots & \vdots \\ 0 & \cdots & g_{nn}\end{pmatrix}g, \mbox{ where } g\in N, \mbox{ and } \nu(g_{ii})=0 \mbox{ for } i=1,\dots,n,
\end{equation*}
we have an element 
\begin{equation*}
g''=\begin{pmatrix}g_{11} & \cdots & 0 \\ \vdots & \ddots & \vdots \\ 0 & \cdots & g_{nn}\end{pmatrix}g\begin{pmatrix}g_{11}^{-1} & \cdots & 0 \\ \vdots & \ddots & \vdots \\ 0 & \cdots & g_{nn}^{-1}\end{pmatrix}
\end{equation*}
in $N$ (here we are using the fact that $N$ is a normal subgroup of $B$) with
\begin{equation*}
dist(g',g'')=l\left( \begin{pmatrix}g_{11}^{-1} & \cdots & 0 \\ \vdots & \ddots & \vdots \\ 0 & \cdots & g_{nn}^{-1}\end{pmatrix}\right) =0.
\end{equation*}

We know that by the virtue of Lemma~\ref{lemmaAsdimN} that $\asdim N=0$, hence $\asdim\tilde N=0$, that is the inverse of the $R$-ball is a finite union of sets of asymptotic dimension $0$. Its asymptotic dimension is $0$ as well, and so is the asymptotic dimension of every $f^{-1}\left( B_R(a)\right)$.

Now we are ready to apply Theorem~\ref{theoremAsdimLipschitz} and conclude that $\asdim B$ does not exceed $\asdim A$, which is $(n-1)$, according to  Lemma~\ref{lemmaAsdimA}, but since $B$ contains $A$ as a subgroup, asymptotic dimension of $B$ is precisely $(n-1)$.
\end{proof}

\subsection{Example: symmetric spaces}\label{sectionAsdimSymmetricSpaces}

The following theorem was proven in~\cite{BellDranishnikov04} by means of connecting the Hirsch length and the asymptotic dimension of a nilpotent group:

\begin{theorem}[{cf.~\cite[Theorem~12]{BellDranishnikov04}}]
Let $G$ be a connected Lie group and let $K$ be its maximal compact subgroup. Then (in the $G$-invariant metric) $\asdim G/K=\dim G/K$.
\end{theorem}

We shall need the straightforward consequence of this theorem:

\begin{corollary}
The symmetric spaces $SL(n,\C)/SU(n)$ and $SL(n,\R)/SO(n)$ have finite asymptotic dimension for any natural $n$.
\end{corollary}

%% file: intchar.tex
\section{Groups of integral characteristic}\label{chapterIntegralCharacteristic}

\begin{definition}
A subgroup $\Gamma$ of $SL(n,\C)$ is said to have an \emph{integral characteristic} if the coefficients of the characteristic polynomial of every element of $\Gamma$ are algebraic integers.
\end{definition}

\begin{remark}
This is equivalent to the condition that the eigenvalues of every element of $\Gamma$ are algebraic integers, since the ring of algebraic integers is integrally closed.
\end{remark}

Typical examples of groups of integral characteristic are unipotent groups, namely subgroups of the conjugates of the upper-triangular matrix groups with $1$'s along the diagonal. However, the entire class of groups of integral characteristic is much richer: for example, it contains all subgroups of $SL(n,\Z)$.

\subsection{Subgroups of $SL(n,\mathbb{Q})$}\label{sectionIntCharSLnQ}

Let $\Gamma$ be a finitely generated subgroup of $SL(n,\mathbb{Q})$, not necessarily of integral characteristic. We think of $\Gamma$ as of a group defined over the ring $A=\mathbb{Z}[\frac1s]$, where $s$ is the l.u.b. of all the denominators of the fractions participating in the entries of the generating set of $\Gamma$. If necessary, we shall enlarge $A$ so that $s$ is the product of distinct primes and therefore we may assume that prime factorization of $s$ is $p_1\cdot p_2\cdots p_m$ with nonrepeating terms.

For any index $k$ between $1$ and $m$ the natural inclusions
\begin{equation*}
\Gamma\subseteq SL(n,A)\subseteq SL(n,\mathbb{Q}_{p_k})
\end{equation*}
allow us to consider the action of $\Gamma$ on the $p_k$-adic building of equivalence classes of lattices in $\mathbb{Q}_{p_k}^n$. Denote this affine building by $X_{p_k}$.

The action $\alpha_{p_k}$ of $SL(n,\Q_{p_k})$ (and, abusing notation, of $\Gamma$ as well) on $X_{p_k}$ is defined via its natural action on lattices.

The action $\alpha_{\infty}$ of $\Gamma$ on the symmetric space $X_{\infty}=SL(n,\R)/SO(n,\R)$ is defined via the natural isometric action of $SL(n,\mathbb{R})$ and the inclusion of $SL(n,A)$ into $SL(n,\R)$.

Finally, $\Gamma$ acts on the product $X$ of buildings and the symmetric space:
\begin{equation*}
X=X_{p_1}\times X_{p_2}\times\dots\times X_{p_m}\times X_{\infty}
\end{equation*}
via the diagonal action
\begin{equation*}
\alpha=\alpha_{p_1}\times \alpha_{p_2}\times\dots\times \alpha_{p_m}\times\alpha_{\infty}
\end{equation*}

\begin{theorem}\label{theoremSLnQ}
The action $\alpha$ is proper.
\end{theorem}

In the proof we shall be using the following: 

\begin{lemma}\label{lemmaBalls}
Let $\Gamma$ be a group acting on a metric space $X$ by isometries. Then the condition
\begin{equation*}
\forall x\in X \qquad \forall C>0 \qquad \#\{g\in\Gamma | \dist(g.x,x)<C\}<\infty
\end{equation*}
implies the condition
\begin{equation*}
\forall x\in X \qquad\qquad \#\{g\in\Gamma | g.B\cap B\ne\emptyset\}<\infty
\end{equation*}
for any ball $B$ in $X$.
\end{lemma}
\begin{proof}
Given ball $B$, we take $x$ to be its center, and $C$ to be $3$ times its radius. Then $B$ and the translated ball $g.B$ do not meet each other, as long as $\dist(g.x,x)\geqslant C$, so that the second set is empty. But if $\dist(g.x,x)<C$, then the first condition guarantees that there are only finitely many group elements with such property.
\end{proof}

\begin{proof}[Proof of the theorem]
We have to show that for any compact $K$
\begin{equation*}
\#\{g\in\Gamma|g.K\cap K\neq\emptyset\}<\infty.
\end{equation*}

Since it is enough to prove the statement for sufficiently large compact sets $K$ only, we shall enlarge the set $K$ in the following way: let $K_{p_k}$ be the projection of $K$ into the $k$-th building, and $K_{\infty}$ be the projection of $K$ into the symmetric space. Clearly
\begin{equation*}
K\subseteq K_{p_1}\times\dots\times K_{p_m}\times K_{\infty},
\end{equation*}
and we shall work with the latter set instead of $K$.

Now the set
\begin{equation*}
\{g\in\Gamma|g.(K_{p_1}\times\dots\times K_{p_m}\times K_{\infty})\cap(K_{p_1}\times\dots\times K_{p_m}\times K_{\infty})\neq\emptyset\}
\end{equation*}
is actually
\begin{equation*}
\bigcap_{k=1}^m\{g\in\Gamma|g.K_{p_k}\cap K_{p_k}\neq\emptyset\}\cap\{g\in\Gamma|g.K_{\infty}\cap K_{\infty}\neq\emptyset\},
\end{equation*}
therefore it is enough to show that the intersection of all $\{g\in\Gamma|g.B\cap B\neq\emptyset\}$ is finite (here $B$ denotes one of $K_{p_1}, \dots, K_{p_m}$, or $K_{\infty}$, and we may assume, by enlarging the compact sets, if necessary, that each $B$ is indeed a closed ball).

Then it suffices to show that for any number $C$ and any points $v_k\in X_{p_k}$ ($k=1,\dots, m$) and $x\in X_{\infty}$
\begin{equation}\label{propernessWithC}
\#\left(\bigcap_{k=1}^m\{g\in\Gamma|\dist_{X_{p_k}}(g.v_k,v_k)<C\}\cap\{g\in\Gamma|\dist_{X_{\infty}}(g.x,x)<C\}\right)<\infty
\end{equation}

Indeed, note that the action of $g$ on all our spaces is isometric, therefore Lemma~\ref{lemmaBalls} can be engaged for each building and the symmetric space, respectively.

By the triangle inequality it is sufficient to check condition \eqref{propernessWithC} for the ``base'' vertices $v_0$ of the buildings and the ``center'' point $x_0$ of the symmetric space.

Now let $g=(g_{ij})$ be an element from the intersection of those sets. Here each matrix entry $g_{ij}$ belongs to $\mathbb{Z}[\frac1{p_1},\dots, \frac1{p_m}]$, namely
\begin{equation*}
g_{ij}=\frac{a_{ij}}{\prod_{k=1}^mp_k^{n_{ijk}}},\qquad\qquad a_{ij}, n_{ijk}\in\mathbb{Z},\qquad (a_{ij},p_k)=1, \qquad k=1,\dots, m.
\end{equation*}

The distance from $g.v_0$ to $v_0$ in the building $X_{p_k}$ is bounded only if each matrix entry has bounded-from-above powers $n_{ijk}$ in the denominator, which means the denominator itself should be bounded.

Since $\dist_{X_{\infty}}(g.x_0,x_0)<C$ means $\cosh(\dist_{X_{\infty}}(g.x_0,x_0))<\cosh C$, and, according to~\cite{BridsonHaefliger},
\begin{equation*}
\cosh(\dist_{X_{\infty}}(g.x_0,x_0))=\sum_{i,j=1}^n g_{ij}^2,
\end{equation*}
each matrix entry should be bounded. Together with the previous observation, this leads to the finite number of choices for $a_{ij}$ and $n_{ijk}$, that is there are finitely many such elements $g$ in the intersection. 
\end{proof}

\subsection{Groups with an irreducible action}\label{sectionIntCharIrrAction}

In the discussion of the next results we fix an algebraically closed field $K$
and a multiplicative monoid $G$ in the matrix algebra $M_n(K)$.

The classical Burnside lemma (see, for example,~\cite{Bass}) claims that:
\begin{lemma}\label{lemmaBurnsideLemma}
If a group $G$ acts irreducibly on $K^n$, then $G$ contains a linear basis $\{g_1, g_2, \dots,
g_{n^2}\}$ of $M_n(K)$.
\end{lemma}

Now we assume that $G$ is a subgroup of $SL(n,K)$ and $\Gamma$ is a subgroup
of $G$ of integral characteristic. We shall construct a matrix representation of $G$ with the property that the image of $\Gamma$ under that representation is arithmetic in most cases. The following result resembles the ``trace'' technique used by Zimmer in his study of Property T groups (see~\cite{Zimmer}), adopted by Higson (unpublished). Denote the field of algebraic numbers by $\A$.

\begin{lemma}\label{lemmaIntCharEmbeddingAlpha}
Let $G$ be a subgroup of $SL(n,K)$. There exists a representation $\alpha:G\to GL(N,\C)$ with $N=n^2$ and such that for every subgroup $\Gamma$ of $G$, which has integral characteristic and acts irreducibly on $K^n$, its image $\alpha(\Gamma)$ is conjugate in $GL(N,\C)$ to a subgroup of $GL(N,\A)$.
\end{lemma}
\begin{proof}
For every $g\in M_n(K)$ define a complex-valued linear functional $f_g$ on
$KM_n(K)=M_n(K)$ by
\begin{equation*}
f_g(h)=\tr(gh).
\end{equation*}
Take some basis $\{g_1, g_2, \dots, g_N\}$ of $M_n(K)$ and consider
\begin{equation*}
V=\Span\{f_{g_1},\dots, f_{g_N}\}.
\end{equation*}
Since the conditions defining $f_g$
are linear with respect to $g$, $\lambda f_g$ can be formally written as
$f_{\lambda g}$, and $f_g+f_h$ as $f_{g+h}$. According to our assumptions $V$
contains every $f_{g_j}$, $j=1, 2, \dots, N$, therefore it contains their
linear span, that is, every $f_g$ with $g\in M_n(K)$. We have $\dim V=N=n^2$ and
\begin{equation*}
V=\Span\{f_{g_1},\dots, f_{g_N}\}=\Span\{f_g\}_{g\in M_n}.
\end{equation*}

Any element $g\in G$ acts on $V$ by
\begin{equation*}
g.f_h=f_{gh}.
\end{equation*}
With respect to the basis $\left\lbrace g_1, g_2, \dots, g_N\right\rbrace $ of $M_n(K)$ and the corresponding basis
$\left\lbrace f_{g_1}, f_{g_2}, \dots, f_{g_N}\right\rbrace $ of $V$ this action is given by a matrix
$(\alpha^g_{ij})$, such that
\begin{equation}\label{matrixElements}
g.f_{g_i}(h)=\sum_{j=1}^N\alpha^g_{ij}f_{g_j}(h)
\end{equation}
Thus we have a representation $\alpha: G\to M_N(\C)$.

To prove the injectivity of $\alpha$ it is enough to show that if
\begin{equation}\label{injectivityCondition}
g.f_1(h)=f_1(h)
\end{equation}
for some $g\in G$ and all $h=g_1,\dots, g_N$, then $g=1$. Taking the linear span over
all such $h$, the condition is equivalent to the one with $h\in M_n(K)$. Plugging
standard elementary matrices for $h$ in \eqref{injectivityCondition}, we
deduce that $g$ is indeed equal to $1$.

In our construction the representation $\alpha$ depends on the choice of the basis  $\left\lbrace g_1,\dots,g_N\right\rbrace $. Suppose $\left\lbrace \tilde g_1,\dots, \tilde g_N\right\rbrace $ is another basis of $M_n(K)$ which gives rise to a representation $\tilde{\alpha}$. Then, since the space $V$ is the same for both bases, and any group element is represented as a linear operator on $V$, the images of any subgroup of $G$ under $\alpha$ and $\tilde{\alpha}$ are conjugate in $GL(N,\C)$. This allows us to choose any convenient basis in the proof of the last statement of the lemma.

Suppose that $\Gamma$ is a subgroup of $G$, has integral characteristic and acts irreducibly on $K^n$. Then the basis $\{g_1, g_2, \dots, g_N\}$ can be taken to be the one furnished by Lemma~\ref{lemmaBurnsideLemma}, that is, consisting of elements of  $\Gamma$.

For $\gamma\in\Gamma$ write \eqref{matrixElements} as
\begin{equation}\label{matrixElementsTr}
\tr(\gamma g_i h)=\sum_{j=1}^N\alpha^{\gamma}_{ij}\tr(g_j h).
\end{equation}

Then $\alpha^{\gamma}_{ij}$ are the solutions of the system of linear
equations with algebraic coefficients, as long as $\gamma\in\Gamma$ (obviously
each $g_j$ is in $\Gamma$ and we can test using $h$ from the basis $\left\lbrace g_1, g_2,
\dots, g_N\right\rbrace $ only), and therefore have to be algebraic as well. This proves the last statement of the lemma.
\end{proof}

\begin{corollary}
Suppose that there exists a finitely generated subring $A$ of $\A$, such that trace of every element in $\Gamma$ lies in $A$. Then $\alpha$ represents $\Gamma$ within $GL(N,\tilde A)$, where $\tilde A$ is a finitely generated subring of the field of fractions of $A$.
\end{corollary}

\begin{proof}
We know that coefficients $\alpha^{\gamma}_{ij}$ are uniquely determined as solutions of \eqref{matrixElementsTr}, and it is enough to check these conditions for $h=g_1, g_2,
\dots, g_N$. Therefore \eqref{matrixElementsTr} is equivalent to a system of linear equations on $\alpha^{\gamma}_{ij}$ with nonzero coefficients accompanying the unknowns  $\tr(g_jg_k)$. Define
\begin{equation*}
\tilde A=A[(\tr(g_1g_1))^{-1}, (\tr(g_1g_2))^{-1}, \dots, (\tr(g_Ng_N))^{-1}].
\end{equation*}
Then \eqref{matrixElementsTr} can be solved by Gauss-Jordan elimination process, involving only ring operations in $\tilde A$, that is all $\alpha^{\gamma}_{ij}$ should belong to $\tilde A$.
\end{proof}

Now we are in a position to prove

\begin{theorem}\label{theoremIntCharIrr}
Given $G$, a finitely generated subgroup of $SL(n,\C)$, there exists an action of $G$ on a finite-asymptotic-dimensional space $X$, such that for every subgroup $\Gamma$ of $G$ of integral characteristic acting irreducibly on $\C^n$ the induced action of $\Gamma$ on $X$ is proper.
\end{theorem}

\begin{proof}
Construct an embedding $\alpha:G\to GL(N,\C)$, as in Lemma~\ref{lemmaIntCharEmbeddingAlpha}. We know that $G$ is finitely generated, thus we can assume that $\alpha:G\to GL(N,F)$, where $F$ is a finitely generated field (generated, say, by all entries of the images of all generators of $G$). The images of elements of any subgroup $\Gamma$ under $\alpha$ are algebraic, and since finitely generated algebraic extension has to be finite, we conclude that $\alpha(\Gamma)\subseteq GL(N,K)$, where $K$ is an extension of $\Q$ of degree $m$.

Now we can identify
\begin{equation*}
GL(N,K) \cong End(K^N) \subseteq End(\Q^{mN}) \cong GL(mN,\Q)
\end{equation*}
and think of $\Gamma$ as of a subgroup of $GL(mN,\Q)$. According to Theorem~\ref{theoremSLnQ} it acts properly on a finite-dimensional space.

We know by the virtue of Lemma~\ref{lemmaIntCharEmbeddingAlpha} that any other subgroup $\tilde\Gamma$ of integral characteristic within $G$ is conjugate to a subgroup of $GL(mN,\Q)$ in our representation, so that it also acts properly.
\end{proof}

\subsection{Diagonal Parts}\label{sectionIntCharDiagonal}

Keep all the notations from Section~\ref{sectionIntCharIrrAction}.

\begin{theorem}
Let $\Gamma$ be a subgroup of $G$ of integral characteristic and such that the restriction of the action $\alpha$ to $\Gamma$ is improper. Then $\Gamma$ is unipotent, that is conjugate to a subgroup of the uni-upper triangular matrices.
\end{theorem}
\begin{proof}
If we have $\Gamma\subseteq G\subset SL(n,K)$, a subgroup of integral characteristic, which does not act irreducibly on $\C^n$, then we can find a tower of linear subspaces of $K^n$
\begin{equation*}
\{0\}=V_0\subset V_1\subset\dots\subset V_r=K^n,
\end{equation*}
such that the induced action of $\Gamma$ on each $V_i/V_{i-1}$ is irreducible, and therefore in a suitable basis every element $g$ of $\Gamma$ is of the form
\begin{equation*}
g=\begin{pmatrix}
g_1 & * & \hdots & * \\
0 &g_2 & \hdots & * \\
\vdots &\vdots & \ddots & \vdots \\
0 &0 & \hdots & g_r \\
\end{pmatrix},
\end{equation*}
where each $g_i$ is the restriction of the natural action of $g$ onto the space $V_i/V_{i-1}$. Consider the map $g\mapsto g_i$. This map is a homomorphism
\begin{equation*}
f_i : \Gamma\to\Gamma_i\subseteq GL(\dim(V_i/V_{i-1}),K),
\end{equation*}
such that $\Gamma_i$ acts irreducibly on $K^{\dim(V_i/V_{i-1})}$, and the results of Section~\ref{sectionIntCharIrrAction} apply to each $\Gamma_i$.

The map
\begin{equation*}
f : \begin{pmatrix}
g_1 & * & \hdots & * \\
0 &g_2 & \hdots & * \\
\vdots &\vdots & \ddots & \vdots \\
0 &0 & \hdots & g_r \\
\end{pmatrix}\mapsto(g_1,g_2,\dots,g_r)
\end{equation*}
is a homomorphism $\Gamma\to\Gamma_1\times\dots\times\Gamma_r$ whose kernel consists of unipotent matrices. Applying the composition of $\alpha$ and $f$ gives us an action which is proper on each $\Gamma_i$, by dint of Lemma~\ref{lemmaIntCharEmbeddingAlpha}.

The isotropy of this action comes from the kernel of the homomorphism $f$, that is it consists of the unipotent matrices.
\end{proof}

We shall proceed to construct a proper action for them in the next section.

\subsection{Unipotent subgroups of $SL(2,\C)$}\label{sectionIntCharUnipotent2}

In the following discussion all fields are supposed to be subfields of $\C$, that is, in particular of characteristic $0$.

Let $K$ be a field, $\theta$ an algebraic over $K$ number of degree $n$, and $\theta_1=\theta,\dots,\theta_n$ be its conjugates. They are all distinct, so that via the well-known formula for the Vandermonde determinant,
\begin{equation*}
\det\begin{pmatrix}
1 & \theta_1 & \dots & \theta_1^{n-1} \\
1 & \theta_2 & \dots & \theta_2^{n-1} \\
\vdots & \vdots & \ddots & \vdots \\
1 & \theta_n & \dots & \theta_n^{n-1} \\
\end{pmatrix}\ne0.
\end{equation*}

This idea can be generalized for the case of polynomials with independent coefficients in the following way (to keep our notation concise, we write $\vec{t}$, possibly with an index, for a finite tuple of indeterminates):

\begin{lemma}\label{lemmaTranscendentalExt}
Let $K$ be a field and $p_1,\dots,p_n\in K[\vec{t}]$ be linearly independent over $K$. Then for any choice of mutually disjoint families $\vec{t}_1,\dots,\vec{t}_n$ of indeterminates
\begin{equation}
\det\begin{pmatrix}
p_1(\vec{t}_1) & \dots & p_n(\vec{t}_1)\\
\vdots & \ddots & \vdots\\
p_1(\vec{t}_n) & \dots & p_n(\vec{t}_n)
\end{pmatrix}\ne0
\end{equation}
\end{lemma}

\begin{proof}
We shall argue by induction on the number $n$ of polynomials. If $n=1$, then the linear independence of $\{p_1\}$ means $p_1\ne0$ and so the determinant in question is nonzero as well.

Now suppose that the statement of the lemma was proven for any choice of $n-1$ linearly independent polynomials. Impose the lexicographic order on the set $\vec{t}$ of indeterminants. Then, working in the $K$-span of $p_1,\dots,p_n$, we can choose a new basis $\tilde{p}_1,\dots,\tilde{p}_n$ for the $K$-span of $p_1,\dots,p_n$, with the property that the leading term of each $\tilde{p}_2,\dots,\tilde{p}_n$ has order strictly less than that of $\tilde{p}_1$. The new system will also be linearly independent over $K$, and the matrix
\begin{equation}
\begin{pmatrix}
\tilde{p}_1(\vec{t}_1) & \dots & \tilde{p}_n(\vec{t}_1)\\
\vdots & \ddots & \vdots\\
\tilde{p}_1(\vec{t}_n) & \dots & \tilde{p}_n(\vec{t}_n)
\end{pmatrix}
\end{equation}
will be obtained from the original matrix by elementary operations on the matrix columns, so that the determinants of two matrices coincide up to a possible sign discrepancy. The determinant of the latter matrix can be decomposed along the first row as
\begin{equation}
\tilde{p}_1(\vec{t}_1)\det\begin{pmatrix}\tilde{p}_2(\vec{t}_2) & \dots & \tilde{p}_n(\vec{t}_2)\\\vdots & \ddots & \vdots\\\tilde{p}_2(\vec{t}_n) & \dots & \tilde{p}_n(\vec{t}_n)\end{pmatrix}+\sum_{i=2}^n(-1)^i\tilde{p}_i(\vec{t}_1)(\mbox{polynomial in }\vec{t}_2,\dots,\vec{t}_n)
\end{equation}

According to the induction hypothesis the coefficient of $\tilde{p}_1(\vec{t}_1)$ of the first term is nonzero, and, since no other term with $\vec{t}_1$ could possibly cancel $\tilde{p}_1(\vec{t}_1)$, the entire determinant is nonzero.
\end{proof}

\begin{remark}\label{remarkTranscendentalExt}
The conclusion of this lemma remains true even if $p_1,\dots,p_n$ lie not in $K[\vec{t}]$, but in the fraction field $K(\vec{t})$. Indeed, there are finitely many denominators in $p_1,\dots,p_n$. After multiplying each $p_i$ by all of them, we obtain $\tilde{p}_1,\dots,\tilde{p}_n$ in $K[\vec{t}]$, linearly independent over $K$ if the original $p_1,\dots,p_n$ were. Now the lemma can be applied to the modified set $\tilde{p}_1,\dots,\tilde{p}_n$, yielding a nonzero determinant. The determinant for $p_1,\dots,p_n$ will then be that nonzero quantity, divided by the product of denominators of $p_1,\dots,p_n$.
\end{remark}

We are ready to start the construction of a proper action. Given a finitely generated group $G\subset SL(2,K)$, where $K$ is a finitely generated subfield of $\C$, we can think of $K$ as of a two-step extension of $\Q$: first we pick some transcendence base $\vec{t}$ of $K$ and obtain a purely transcendental extension $\Q(\vec{t})$, and then $K$ itself is an algebraic extension of $\Q(\vec{t})$, moreover, a finite extension, and so it is generated by an element $\theta$, algebraic over $\Q(\vec{t})$ of degree $n$.

Given such $K$, fix a natural number $m$ which will serve as a common bound on the dimensions of unipotent subgroups of $G$ and pick $mn\#\{\vec{t}\}$ complex numbers, algebraically independent over $K$. We shall group these numbers into $mn$ tuples of $\#\{\vec{t}\}$ elements each and denote those tuples as
\begin{equation*}
\vec{t}_1, \vec{t}_2, \dots, \vec{t}_{mn},
\end{equation*}
and so set up a component-wise bijection between each $\vec{t_j}$ and the transcendence base $\vec{t}$. The existence of these numbers follows from the following inductive procedure: take the algebraic closure of $K$; it is countable since $K$ is. Now pick the first number from the complement of the closure and extend the closure by adding this element. The new field is again countable, and so is its closure, and we can pick the second complex number from its complement, and so on.

Let $\theta=\theta_1, \theta_2, \dots, \theta_n$ be conjugates of $\theta$. Define $n$ embeddings
\begin{equation*}
\sigma_i: K=\Q(\vec{t})(\theta)\to\Q(\vec{t})(\theta_i)\subset\C, \qquad\qquad i=1,\dots,n
\end{equation*}
by letting $\sigma_i$ to be an identity on $\Q(\vec{t})$ and sending $\theta$ to $\theta_i$.

Also define $mn$ embeddings
\begin{equation*}
\sigma_{i+n}: K=\Q(\vec{t})\to\Q(\vec{t_i})\subset\C, \qquad\qquad i=1,\dots,mn
\end{equation*}
by identifying $\vec{t}$ with $\vec{t_i}$ component-wise. Extend these embeddings to $K$ by mapping $\theta$ to some number, algebraic over $\Q(\vec{t_i})$ of degree $n$.

Now let $X\times X\times\dots\times X$ be a product of $mn+n$ copies of the symmetric space $SL(2,\C)/SU(2)$ and define an action of $SL(2,K)$ on this space by
\begin{equation}\label{equationTwistedAction2}
g.(x_1,x_2,\dots,x_{mn+n})=(\sigma_1(g).x_1,\sigma_2(g).x_2,\dots,\sigma_{mn+n}(g).x_{mn+n}),
\end{equation}
where $\sigma_i(g).x_i$ on the right-hand side means applying $\sigma_i$ to every entry of $g\in SL(2,K)$ and making the resulting matrix in $SL(2,\C)$ act on the point $x_i$ of the $i$-th symmetric space in the usual way.

\begin{theorem}\label{theoremIntCharUnipotent}
The restriction of the action~\eqref{equationTwistedAction2} to any unipotent subgroup of $G$ of dimension $m$ is proper.
\end{theorem}

\begin{proof}
Any unipotent subgroup of dimension $m$ is conjugate to
\begin{equation*}
\begin{pmatrix}1 & u_1\Z+\dots+u_m\Z \\ 0 & 1 \end{pmatrix},
\end{equation*}
where $u_1,\dots,u_m$ are linearly independent over $\Q$ elements of $K$. Let
\begin{equation*}
u_j=\sum_{l=0}^{n-1}p_{j,l}(\vec{t})\theta^l, \qquad\qquad p_{j,l}(\vec{t})\in \Q(\vec{t}), j=1,\dots,m.
\end{equation*}

Seeking a contradiction, assume that the restriction of the action~\eqref{equationTwistedAction2} is improper. Then there exist infinitely many different $m$-tuples of integers $(z_1,\dots,z_m)$, such that the quantities
\begin{equation}\label{equationSigmaUs}
\sigma_i(z_1u_1+\dots+z_mu_m), \qquad\qquad i=1,\dots,mn+n
\end{equation}
are all bounded. For $i=1,\dots,n$ write this condition as
\begin{equation*}
\sum_{j=1}^mz_j\sum_{l=0}^{n-1}p_{j,l}(\vec{t})\theta_i^l=\mbox{bounded}, \qquad\qquad i=1,\dots,n,
\end{equation*}
that is
\begin{equation*}
\begin{pmatrix}
1 & \theta_1 & \dots & \theta_1^{n-1} \\
1 & \theta_2 & \dots & \theta_2^{n-1} \\
\vdots & \vdots & \ddots & \vdots \\
1 & \theta_n & \dots & \theta_n^{n-1} \\
\end{pmatrix}\begin{pmatrix}
\sum_{j=1}^mz_jp_{j,0}(\vec{t}) \\
\sum_{j=1}^mz_jp_{j,1}(\vec{t}) \\
\vdots \\
\sum_{j=1}^mz_jp_{j,n-1}(\vec{t}) \\
\end{pmatrix}=\mbox{ bounded}.
\end{equation*}

As we pointed out before, the $\theta$-matrix has a (fixed) inverse, so that the sum $\sum_{j=1}^mz_jp_{j,l}(\vec{t})$ must be bounded for every $l=0,\dots,n-1$.

Write this condition as
\begin{equation*}
\begin{pmatrix}
p_{1,0}(\vec{t}) & \cdots & p_{m,0}(\vec{t}) \\
\vdots & \ddots & \vdots \\
p_{1,n-1}(\vec{t}) & \cdots & p_{m,n-1}(\vec{t}) \\
\end{pmatrix}
\begin{pmatrix}
z_1 \\
\vdots \\
z_m
\end{pmatrix}=
P(\vec{t})\begin{pmatrix}
z_1 \\
\vdots \\
z_m
\end{pmatrix}=\mbox{ bounded}.
\end{equation*}

Transform the matrix $P(\vec{t})$ into its echelon form $\tilde P(\vec{t})$ by means of elementary operations with the rows over $\Q(\vec{t})$. Note that the linear independence of $u_1,\dots,u_m$ over $\Q$ means linear independence of columns of such a matrix, even in echelon form. In other words, the obtained echelon form has at least one nonzero row.

It could happen that some nonzero entries of the resulting matrix are rational. In this case, multiply the rows where this happens by some elements of $\Q(\vec{t})$ to ensure that all nonzero entries are in $\Q(\vec{t})\backslash\Q$ and continue to denote the resulting matrix by $\tilde P(\vec{t})$.

Now apply $\sigma_{i+n}$ to the $i$-th row of $\tilde P(\vec{t})$ (that is, write $\vec{t}_i$ instead of $\vec{t}$ in that row) and add all the rows together to obtain a row which we denote as
\begin{equation}\label{equationAddEmbeddings}
(\tilde p_1(\vec{t}_1,\dots,\vec{t}_n),\dots,\tilde p_m(\vec{t}_1,\dots,\vec{t}_n)).
\end{equation}

We claim that the elements $\tilde p_1,\dots,\tilde p_m$ are linearly independent over $\Q$. Indeed, suppose that there exist $m$ rational numbers $\lambda_1,\dots,\lambda_m$, not all of them zeros, such that $\lambda_1\tilde p_1(\vec{t}_1,\dots,\vec{t}_n)+\dots+\lambda_m\tilde p_m(\vec{t}_1,\dots,\vec{t}_n)=0$. Write this condition as
\begin{equation*}
\sum_{j=1}^m\lambda_j\sum_{i=1}^n\tilde p_{j,i-1}(\vec{t}_i)=0.
\end{equation*}
Note that each $\tilde p_{j,i-1}(\vec{t}_i)$ is either $0$ or algebraically independent from other $\tilde p_{j,i-1}$'s with the same $i$. This means
\begin{equation*}
\sum_{j=1}^m\lambda_j\tilde p_{j,i-1}(\vec{t}_i)=0
\end{equation*}
for all indices $i$ for which not all $\tilde p_{j,i-1}(\vec{t}_i)$ are zeros. This implies that the columns of $\tilde P$ and $P$ as well are linearly dependent with the same choice of coefficients $\lambda_1,\dots,\lambda_m$, contradictory to our assumption.

Now we apply $\sigma_{i+2n}, \sigma_{i+2n}, \dots, \sigma_{i+mn-n}$ to the $i$-th row of $\tilde P(\vec{t})$ and add rows together. As a result, we obtain $m$ rows, similar to~\eqref{equationAddEmbeddings}, with the only difference that the variables in the $j$-th row are $\vec{t}_{n(j-1)+1},\dots,\vec{t}_{nj}$. Thus we have:
\begin{equation}\label{equationBoundness}
\begin{pmatrix}
p_1(\vec{t}_1,\dots,\vec{t}_n) & \cdots & p_m(\vec{t}_1,\dots,\vec{t}_n) \\
\vdots & \ddots & \vdots \\
p_1(\vec{t}_{mn-n+1},\dots,\vec{t}_{mn}) & \cdots & p_m(\vec{t}_{mn-n+1},\dots,\vec{t}_{mn}) \\
\end{pmatrix}\begin{pmatrix}
z_1 \\
\vdots \\
z_m
\end{pmatrix}=\mbox{ bounded}.
\end{equation}

Lemma~\ref{lemmaTranscendentalExt} says that the matrix on the left-hand side of~\eqref{equationBoundness} has a fixed inverse, which means $z_1,\dots,z_m$ are all bounded, yielding a contradiction.
\end{proof}

\begin{remark}\label{remarkIntCharUnipotentBounded}
Let us summarize that in the proof of the previous theorem we have shown that if for any unipotent subgroup the $(1,2)$ matrix entry of every element of it contains only $m$ linearly independent over $\Q$ numbers, then the boundness of this matrix entry under all embeddings implies that this matrix entry actually takes a finite number of possible values.
\end{remark}

This discussion culminates in
\begin{theorem}\label{theoremIntCharUnipotent2}
If $\Gamma$ is a finitely generated subgroup of $SL(2,\C)$, such that there exists a uniform bound on $\asdim$ of all unipotent subgroups of $\Gamma$, then there is a finite-asymptotic-dimensional space, on which $\Gamma$ acts with the restriction of the action on any unipotent subgroup of $\Gamma$ being proper.
\end{theorem}

\begin{proof}
Let $m$ be the common bound on $\asdim$ of all unipotent subgroups of $\Gamma$, and $X\times\dots\times X$ be the product of symmetric spaces, equipped with the action~\eqref{equationTwistedAction2}. Then for any unipotent subgroup of $\Gamma$ of dimension not more than $m$ we have enough embeddings to apply Theorem~\ref{theoremIntCharUnipotent}.
\end{proof}

\subsection{Unipotent subgroups of $SL(n,\C)$}\label{sectionIntCharUnipotent}

Suppose we have $G$, a unipotent subgroup of some finitely generated subgroup of $SL(n,K)$. We shall reduce the construction of a proper action for this group to the case of $SL(2,K)$. Assume for now that $G$ is uni-upper-triangular:
\begin{equation}\label{equationCanonicalUnipotent}
G\subseteq\{[g_{ij}] \; | \; g_{ii}=1, g_{ij}=0 \mbox{ for all } i=1,\dots,n, \quad j<i \}.
\end{equation}

We start by splitting $G$ into separate ``layers'' $G_0=G, G_1,\dots,G_{n-2}$ in the inductive fashion which we shall describe in a moment, but first we give an informal illustration of this process for the case $n=4$. Our goal is to get the following sequence:
\begin{equation*}
G_0=G=\begin{pmatrix}1&*&*&*\\0&1&*&*\\0&0&1&*\\0&0&0&1\end{pmatrix}\triangleright
G_1=\begin{pmatrix}1&0&*&*\\0&1&0&*\\0&0&1&0\\0&0&0&1\end{pmatrix}\triangleright
G_2=\begin{pmatrix}1&0&0&*\\0&1&0&0\\0&0&1&0\\0&0&0&1\end{pmatrix}.
\end{equation*}
(The subscript of $G$ corresponds to the thickness of the band with zeros above the diagonal with units.)

Now we give formal definitions. Let $G_1$ be the normal subgroup of $G$ comprised of matrices $[g_{ij}]$ for which $g_{ij}=0$ for $j=i+1$. 

Define a group homomorphism
\begin{gather*}
\phi_1:G\to\prod_{i=1}^{n-1}\left\{\left.\m{1}{x}{0}{1} \right| x\in K \right\}=H_1,\\
\phi_1([g_{ij}])=\left(\m{1}{g_{12}}{0}{1},\dots,\m{1}{g_{n-1,n}}{0}{1}\right).
\end{gather*}

Suppose that we have already defined $G_1,\dots,G_{k-1}$ for some $k\leqslant n-2$. Let $G_k$ be the (normal) subgroup of $G_{k-1}$ of the following kind:
\begin{equation*}
G_k=\left\{[g_{ij}]\in G_{k-1} \quad | \quad g_{ij}=0 \mbox{ for } j=i+k\right\}.
\end{equation*}

Also define a group homomorphism
\begin{gather*}
\phi_k:G_{k-1}\to\prod_{i=1}^{n-k}\left\{\left.\m{1}{x}{0}{1} \right| x\in K \right\}=H_k,\\
\phi_k([g_{ij}])=\left(\m{1}{g_{1,1+k}}{0}{1},\dots,\m{1}{g_{n-k,n}}{0}{1}\right).
\end{gather*}

\begin{definition}\label{definitionBoundedlyComposed}
Let $m$ be a natural number. We say that a uni-upper-triangular group $G$ has \emph{composition rank $m$} if in the notation introduced above,
\begin{enumerate}
\item $\asdim G_{n-2}\leqslant m$.\label{itemBClast}
\item $\asdim \phi_k(G_{k-1})\leqslant m$ for $k=1,\dots,n-2$.\label{itemBCphi}
\end{enumerate}
\end{definition}

For a general unipotent group $G$, we define its composition rank to be the rank of the corresponding uni-upper-triangular group (the one our group is conjugate to).

The idea of this definition is to introduce an effective way to test finite asymptotic dimensionality of $G$. The following paragraphs describe the structure of a group with bounded composition rank and also introduce some notation which we shall use later.

Start with $G_{n-2}$. Since it has asymptotic dimension at most $m$, there are $m$ elements of $K$, which we shall call $u_1^{(1,n)},\dots,u_m^{(1,n)}$, such that
\begin{equation*}
G_{n-2}=\left\{\begin{pmatrix}
1 & 0 & \cdots & 0 & u_1^{(1,n)}\Z+\dots+u_m^{(1,n)}\Z \\
0 & 1 & \cdots & 0 & 0 \\
\vdots & \vdots & \ddots & \vdots & \vdots \\
0 & 0 & \cdots & 0 & 1
\end{pmatrix}\right\}.
\end{equation*}
(more precisely, we should say that there are not more than $m$ linearly independent over $\Q$ elements of $K$ with such representation, but we can allow some ``fake'' generators $u_*^{(*,*)}$ to simplify the notation.)

At the same time, the image $\phi_{n-2}$ of $G_{n-3}$ in $H_2$ has asymptotic dimension not more than $m$ by definition, which means that there exist two $m$-tuples
\begin{equation*}
u_1^{(1,n-1)},\dots,u_m^{(1,n-1)} \quad\mbox{ and }\quad u_1^{(2,n)},\dots,u_m^{(2,n)}
\end{equation*}
of elements of $K$ (again, we let them have the same cardinality $m$ without requiring linear independence), such that $\phi_{n-2}(G_{n-3})$ lies in
\begin{equation*}
\m{1}{u_1^{(1,n-1)}\Z+\dots+u_m^{(1,n-1)}\Z}{0}{1}\times\m{1}{u_1^{(2,n)}\Z+\dots+u_m^{(2,n)}\Z}{0}{1}.
\end{equation*}

Treating each ``layer'' $G_k$ of $G$ in this fashion, we obtain $\sum_{i=1}^{n-1}i=\frac{n(n-1)}2$ $m$-tuples
\begin{equation}\label{equationTuplesOfGenerators}
\{u_1^{(i,j)},\dots,u_m^{(i,j)}\}, \qquad\qquad i=1,\dots,n-1, \qquad j=n-i+1,\dots,n,
\end{equation}
such that for $j=2,\dots,n$ the image $\phi_{j-1}(G_{j-2})$ is in
\begin{equation*}
\m{1}{u_1^{(1,j)}\Z+\dots+u_m^{(1,j)}\Z}{0}{1}\times\dots\times\m{1}{u_1^{(n-j+1,n)}\Z+\dots+u_m^{(n-j+1,n)}\Z}{0}{1}.
\end{equation*}

Now we are ready to describe the matrix entries of elements of $G$. 

\begin{lemma}\label{lemmaStructureUniUpperTriangular}
A uni-upper-triangular group $G$ with composition rank $m$ is encoded by the $m$-tuples~\eqref{equationTuplesOfGenerators} in the following sense. For any group element $[g_{ij}]\in G$ its matrix entries satisfy the following conditions:
\begin{gather*}
g_{ij}=0 \qquad \mbox{ for } i>j, \\
g_{ij}=1 \qquad \mbox{ for } i=j, \\
g_{ij}\in u_1^{(i,j)}\Z+\dots+u_m^{(i,j)}\Z+P_{ij}\qquad \mbox{ for } i<j,
\end{gather*}
where
\begin{equation*}
P_{ij}=\sum\left\{\left.u_k^{(\tilde i,\tilde j)}u_l^{(\hat i,\hat j)}\Z \quad \right| \quad k,l=1,\dots,m,\quad \tilde j-\tilde i<j-i>\hat j-\hat i\right\}.
\end{equation*}
\end{lemma}

\begin{proof}
The first two conditions are trivial.

For the last one note that the condition describing $P_{ij}$ means all possible $\Z$-combinations of the pairwise products of $u_*^{(*,*)}$ taken from $\phi_{\tilde k}(G_{\tilde k-1})$ for all $\tilde k$ less than $k$ mapping $g_{ij}$ into its spot in $H_k$. Now proceed by induction. For elements $g_{ij}$ with $j=i+1$ (that is lying along the first diagonal above the $1$'s) we can only pick a $\Z$-combination of the corresponding $u_*^{(i,j)}$, due to the nature of matrix multiplication. For the next diagonal $g_{ij}$ can still get elements $u_*^{(i,j)}$, plus the portion coming from the \emph{previous} diagonal (that's what $P_{ij}$ stands for). For the diagonal after that one we can get entries from $u_*^{(i,j)}$ and also from two previous diagonals, and so on.
\end{proof}

It is natural to ask how the properties of having bounded composition rank and having finite asymptotic dimension are connected. We appeal to the following result of Bell and Dranishnikov.

\begin{theorem}[{\cite[Theorem~7]{BellDranishnikov04}}]\label{theoremBDSurjection}
Let $\phi:G\to H$ be a surjective homomorphism of a finitely generated group with kernel $K$. Then $\asdim G\leqslant \asdim H+\asdim K$.
\end{theorem}

\begin{corollary}\label{corollaryBoundedlyComposedFAD}
Groups with bounded composition rank have finite asymptotic dimension.
\end{corollary}

\begin{proof}
Keep all the notation used in Definition~\ref{definitionBoundedlyComposed}. We shall prove that $\asdim G_k<\infty$ inductively for $k=n-2,\dots,0$.

Condition~\ref{itemBClast} of Definition~\ref{definitionBoundedlyComposed} means that $\asdim G_{n-2}<\infty$.

For the inductive step suppose that $\asdim G_k<\infty$ has been proven for some $k$. Let $H$ be the image of $G_{k-1}$ under $\phi_k$. According to condition~\ref{itemBCphi} of Definition~\ref{definitionBoundedlyComposed}, $\asdim H<\infty$. Applying Theorem~\ref{theoremBDSurjection} for homomorphism $\phi_k$, its image $H$ and kernel $G_k$ we conclude that $\asdim G_{k-1}<\infty$.

Since $G_0=G$, we are done.
\end{proof}

Finally we are in a position to present the construction of a proper action. We shall employ the same field embeddings and the twisted action, as in Section~\ref{sectionIntCharUnipotent2}. For convenience of the reader a brief reminder of that construction follows, with a few cosmetic changes.

We are given a finitely generated subfield $K$ of $\C$ and the uniform bound $m$ on the composition rank of all unipotent subgroups of a given finitely generated subgroup of $SL(n,K)$.

Let $\vec{t}$ be a transcendence base of $K$ over $\Q$, so that $K$ is an algebraic extension of degree $r$ of $\Q(\vec{t})$. Pick $mr\#\{\vec{t}\}$ algebraically independent over $K$ complex numbers and group them into $mr$ tuples, each of $\#\{\vec{t}\}$ elements:
\begin{equation}\label{equationBasesTuples}
\vec{t}_1,\dots,\vec{t}_{mr}.
\end{equation}

Define $r$ embeddings of $K$ into $\C$, fixing $\Q(\vec{t})$ and mapping the generator of $K$ over $\Q(\vec{t})$ to its conjugates in this extension. Call these embeddings $\sigma_1,\dots,\sigma_r$.

Define $mr$ embeddings $\Q(\vec{t})$ into $\C$ by fixing $\Q$ and sending $\vec{t}$ to $\vec{t}_i$ from the list~\eqref{equationBasesTuples}. Extend each of them to $K$ by sending the algebraic generator of $K$ over $\Q(\vec{t})$ to some complex number of degree $r$ over the corresponding $\Q(\vec{t}_i)$. Call the extensions $\sigma_{r+1},\dots,\sigma_{mr+r}$.

Finally, let $X\times X\times\dots\times X$ be a product of $mr+r$ copies of the symmetric space $SL(n,\C)/SU(n)$, each equipped with a natural action of $SL(n,\C)$. Define a twisted action of $SL(n,K)$ on the product of symmetric spaces in the following way:
\begin{equation}\label{equationTwistedAction}
g.(x_1,x_2,\dots,x_{mr+r})=(\sigma_1(g).x_1,\dots,\sigma_{mr+r}(g).x_{mr+r}), \quad g\in SL(n,K).
\end{equation}

\begin{theorem}[cf. Theorem~\ref{theoremIntCharUnipotent}]\label{theoremIntCharUnipotentProper}
The restriction of the action~\eqref{equationTwistedAction} to any unipotent subgroup of our given subgroup of $SL(n,K)$ is proper.
\end{theorem}

\begin{proof}
As in the proof of Theorem~\ref{theoremIntCharUnipotent} we assume that the action is improper and deduce that this implies that there are infinitely many unipotent group elements with uniformly bounded under all embeddings $\sigma_1,\dots,\sigma_{mr+r}$ matrix entries.

Every unipotent subgroup of our group is conjugate to some uni-upper-triangular subgroup, and since such a conjugation is given by a fixed matrix, we have infinitely many uni-upper-triangular group elements with entries, uniformly bounded under all embeddings.

Keep the notations from Lemma~\ref{lemmaStructureUniUpperTriangular}. The boundness of the $(i,j)$-th entry for $j-i=1$ in the light of Remark~\ref{remarkIntCharUnipotentBounded} means that this entry is expressible as
\begin{equation*}
z_1u_1^{(i,j)}+\dots+z_mu_m^{(i,j)}, \quad (z_1,\dots,z_m)\in\mbox{bounded subset of }\Z^m.
\end{equation*}
This means that we can expect only finitely many different numbers for this line $j-i=1$.

Consider the entries along the next line, $j-i=2$. We know that each of them is bounded. In addition, the $P_{ij}$ portion of each entry comes from the previous diagonal, in our case $j-i=1$, and, since we have already proven that there are finitely many options there, the set of possible offsets $P_{ij}$ is bounded. Hence the linear part
\begin{equation*}
z_1u_1^{(i,j)}+\dots+z_mu_m^{(i,j)}, \quad (z_1,\dots,z_m)\in\Z^m
\end{equation*}
is bounded as well, that is, invoking Remark~\ref{remarkIntCharUnipotentBounded} again, the set of all possible coefficients $(z_1,\dots,z_m)$ is bounded for every entry along this line. Summarizing, only finitely many different entries along the line $j-i=2$ are possible.

Inductively we show that there are only finitely many options for each matrix entry, but this contradicts the assumption that we start with an infinite subset of group elements.
\end{proof}

%% file: maintheorem.tex
\section{Proof of the main theorem}\label{chapterMainTheorem}

Let $\Gamma$ be a finitely generated subgroup of $SL(n,\C)$. Since $\Gamma$ is finitely generated, we can think of it as of a subgroup of $SL(n,A)$, where $A$ is a finitely generated subring of $\C$. Let $K$ be the quotient field of $A$.

The following construction was introduced by Alperin and Shalen in their study of the cohomological dimension on linear groups (see~\cite{AlperinShalen}). We start by recalling some basic facts about integral ring extensions. Within the scope of this discussion all rings under consideration are subrings of $\C$.

The following result is known as the Noether Normalization Lemma.

\begin{lemma}\label{lemmaNoetherNormalizationLemma}
Let $K$ be a field and $A$ an integral domain, finitely generated by the elements $y_1,\dots,y_n$ as an algebra over $K$. Then there are elements $x_1,\dots,x_r\in A$, algebraically independent over $K$, such that $A$ is integral over $K[x_1,\dots,x_r]$.
\end{lemma}

Note that the number $r$ in Lemma~\ref{lemmaNoetherNormalizationLemma} corresponds to the transcendence degree of $A$ over $K$. We shall use this lemma in the form of

\begin{corollary}\label{corollaryNoetherNormalization}
Let $A$ be a finitely generated integral domain. Then there exists a nonzero integer $s$ and algebraically independent over $\Q$ elements $x_1,\dots,x_r\in A[s^{-1}]$, such that $A[s^{-1}]$ is an integral extension of $\Z[s^{-1}][x_1,\dots,x_r]$.
\end{corollary}

\begin{proof}
Let $y_1,\dots,y_n$ be the generating set for $A$ and let $B$ be the ring generated by $y_1,\dots,y_n$ over $\Q$. Apply Lemma~\ref{lemmaNoetherNormalizationLemma} for the ring $B$ and the field $\Q$ and obtain  $x_1,\dots,x_r\in B$ with the property that $B$ is an integral extension of $\Q[x_1,\dots,x_r]$.

Each $y_i$, as an element of $B$, is a root of a monic polynomial with coefficients in $\Q[x_1,\dots,x_r]$. Take $s_1\in\Z$ to be the common denominator of all rationals, participating in the coefficients of the corresponding polynomials for all $y_i$'s. Then each $y_i$ is integral over $\Z[s_1^{-1}][x_1,\dots,x_r]$.

It could happen that not all $x_j$'s lie in $A$. We know that each of them is a $\Q$-linear combination of the products of $y_i$'s, so let $s_2\in\Z$ be the common denominator of all rationals participating in such expressions. Then each $x_j$ lies in the ring extension of $\Z[s_2^{-1}]$ by $y_1,\dots,y_n$, which is $A[s_2^{-1}]$.

Finally, let $s=s_1s_2$. Then each generator $y_i$ of $A$ is integral over $\Z[s^{-1}][x_1,\dots,x_r]$, moreover $\Z[s^{-1}][x_1,\dots,x_r]\subseteq A[s^{-1}]$, and so the latter is an integral extension of the former.
\end{proof}

\begin{lemma}\label{lemmaMainThmValuations}(``Alperin-Shalen lemma'')
Let $A$ be a finitely generated ring with field of fractions $K$. Then there exist finitely many discrete valuations $\nu_1,\dots,\nu_m$ on $A$ with
$$A\cap\mathcal O_{\nu_1}\cap\dots\cap\mathcal O_{\nu_m}\subset\mathcal O,$$
where $\mathcal O$ is the ring of algebraic integers  in $K$.
\end{lemma}

\begin{proof}
We start by finding $s\in\Z$ and $x_1,\dots,x_r$ in $A[s^{-1}]$ such that the polynomial ring $\left(\Z[s^{-1}]\right)[x_1,\dots,x_r]$ has an integral extension $A[s^{-1}]$. For each $i=1,\dots,r$ the ring $\Z[s^{-1}][x_1,\dots,x_r]$, as an extension of $\Z[s^{-1}][x_1,\dots,x_{i-1},x_{i+1},\dots,x_r]$, admits a discrete valuation $\mu_i$, given by
\begin{equation*}
\mu_i(f)=-\deg f, \qquad f(x)\in\Z[s^{-1}][x_1,\dots,x_{i-1},x_{i+1},\dots,x_r][x].
\end{equation*}

These valuations can be extended to $A[s^{-1}]$ and also to $K$ in finitely many ways (see, for example,~\cite[Chapters 3 and 4]{McCarthy}) to produce discrete valuations
$\nu_1,\dots,\nu_t$.

If $x\in A[s^{-1}]\cap\mathcal O_{\nu_1}\cap\dots\cap\mathcal O_{\nu_t}$, its minimal polynomial over $\Q(x_1,\dots,x_r)$ has coefficients in $\left(\Z[s^{-1}]\right)[x_1,\dots,x_r]\cap\mathcal O_{\mu_1}\cap\dots\cap\mathcal O_{\mu_r}$. Thus $x$ is integral over $\Z[s^{-1}]$.

Finally add $p$-adic valuations $\nu_{t+1},\dots,\nu_m$, which correspond to all prime divisors of $s$.

Now if $x\in A\cap\mathcal O_{\nu_1}\cap\dots\cap\mathcal O_{\nu_m}$, it must be integral over $\Z$ and therefore $x\in\mathcal O$.
\end{proof}

\begin{theorem}\label{theoremMainThmAS}
For any finitely generated group $\Gamma\subset SL(n,A)$ there exists a finite product of $(n-1)$-dimensional affine buildings, each of finite asymptotic dimension, on which $\Gamma$ acts in such a way that the isotropy groups lie among the subgroups of $\Gamma$ of integral characteristic.
\end{theorem}

\begin{proof}
Starting from a finitely generated $\Gamma\subset SL(n,A)$, one can find the discrete valuations $\nu_1,\dots,\nu_m$ on $K$, described in Lemma~\ref{lemmaMainThmValuations} above. Let
\begin{equation*}
X=X_{\nu_1}\times X_{\nu_2}\times\dots\times X_{\nu_m},
\end{equation*}
where $X_{\nu_i}$ is the $(n-1)$-dimensional affine building, corresponding to a discrete  valuation $\nu_i$.

The entire group $SL(n,K)$ acts on each $X_{\nu_i}$ in the obvious way, via the natural action on the lattices in $K^n$. Hence $\Gamma$ acts on $X_{\nu_i}$ as well. Finally, define the diagonal action of $\Gamma$ on $X$:
\begin{equation}\label{equationActionAS}
g.(x_1, x_2, \dots, x_m)=(g.x_1, g.x_2, \dots, g.x_m), \qquad\qquad g\in\Gamma, x_i\in X_{\nu_i}.
\end{equation}

The isotropy subgroups of this action are the subgroups of $\Gamma$ which stabilize a vertex in each $X_{\nu_i}$. We know that the stabilizer of any vertex in $X_{\nu_i}$ stabilizes the associated class of $\mathcal O_{\nu_i}$-lattices in $K^n$, thus the coefficients of its characteristic polynomial are in $\mathcal O_{\nu_i}$. Therefore the isotropy elements of the action~\eqref{equationActionAS} have the coefficients of their characteristic polynomials in $\mathcal O_{\nu_1}\cap\dots\cap\mathcal O_{\nu_m}$, and according to Lemma~\ref{lemmaMainThmValuations} they are algebraic integers in $A$, which means that the isotropy groups have integral characteristic.
\end{proof}

Finally we are able to prove Theorem~\ref{theoremMain}:


\begin{proof}[Proof of the main theorem]
Take $X$ to be the product of the buildings furnished in the proof of Theorem~\ref{theoremMainThmAS} and the finite-asymptotic-dimensional spaces from Theorem~\ref{theoremIntCharIrr} and Theorem~\ref{theoremIntCharUnipotent2}. Then the diagonal action of $\Gamma$ on $X$ is proper. Indeed, the stabilizer of such action should be of integral characteristic (Theorem~\ref{theoremMainThmAS}), without any irreducible diagonal parts (Theorem~\ref{theoremIntCharIrr}), that is unipotent. But all such subgroups act properly because of Theorem~\ref{theoremIntCharUnipotent2}.
\end{proof}

%% file: sl2c.tex
\section{Digression to the case of $SL(2)$}

In this final section we shall give a more detailed treatment of linear groups of $2\times2$ matrices by pushing the techniques developed in the previous sections further.

In what follows, we shall extensively use the fact that the Baum-Connes conjecture (with coefficients) passes to subgroups, without explicitly mentioning it.

\subsection{Subgroups of $SL(2,\Q)$}\label{SubgroupsOfSL2Q}

Let $\Gamma$ be a finitely generated subgroup of $SL(2,\Q)$. As before, $\Gamma$ is defined over the ring $A=\Z[\frac1s]$, where $s$ is the l.u.b. of all the denominators of the fractions participating in the entries of the generating set. We may assume that prime factorization of $s$ is $p_1\cdot p_2\cdots p_n$ with nonrepeating terms.

As in Section~\ref{sectionIntCharSLnQ}, consider the trees ($1$-dimensional buildings)  $T_{p_k}$ corresponding to each $p_k$-adic valuation on $A$ and denote by $\alpha_{p_k}$ the induced actions of $\Gamma$ on these trees.

Define an action $\alpha_H$ of $\Gamma$ on the 2-dimensional real hyperbolic space $\HH_2$ via the natural isometric action of $SL(2,\R)$ on $\HH_2$.

Finally, $\Gamma$ acts on the product of trees and a hyperbolic space $\HH_2$:
\begin{equation*}
T=T_{p_1}\times T_{p_2}\times\dots\times T_{p_n}\times\HH_2
\end{equation*}
via the diagonal action
\begin{equation*}
\alpha=\alpha_{p_1}\times \alpha_{p_2}\times\dots\times \alpha_{p_n}\times\alpha_H.
\end{equation*}

\begin{theorem}\label{ProperActionSL2Q}\label{theoremProperActionSL2Q}
The action $\alpha$ is proper.
\end{theorem}

\begin{proof}
The proof of Theorem~\ref{theoremSLnQ} is still applicable here, with a minor change: instead of a symmetric space $SL(n,\C)/SU(n)$ we have a hyperbolic space $\HH_2$.

Since for the hyperbolic space and its ``center'' $x_0$ the property
\begin{equation*}
\cosh(\dist_{\mathbb{H}_2}(g.x_0,x_0))=\sum_{i,j=1}^2 g_{ij}^2,\qquad g\in SL(n,\R)
\end{equation*}
still holds, the rest of the argument from the proof of Theorem~\ref{theoremSLnQ} works without any alteration.
\end{proof}

\subsection{Reduction to groups of integral characteristic}\label{ReductionToIntegral}

In this section we construct a hierarchy of families of subgroups of $SL(2,\C)$ in such a way that its base consists of groups of integral characteristic, each subgroup from some level of the hierarchy admits an action on a tree with isotropy belonging to the previous level, and show that 
any finitely generated subgroup of $SL(2,\C)$ belongs to some level of this hierarchy.

The main motivation for the study of the isotropy of the group action on a tree is the following

\begin{theorem}[Oyono-Oyono,~\cite{OyonoOyono}]\label{OyonoOyonoTheorem}\label{theoremOyonoOyono}
Let $\Gamma$ be a discrete countable group acting on a tree. Then the Baum-Connes conjecture holds for $\Gamma$ if and only of it holds for all the isotropy subgroups of the action on vertices of the tree.
\end{theorem}

Define a sequence
$$\mathcal H_0,\mathcal H_1,\dots,\mathcal H_i,\dots$$
of families of subgroups of $SL(2,\mathbb C)$ in the following way:
\begin{itemize}
\item $\mathcal H_0$ consists of groups of integral characteristic.
\item $\mathcal H_i$ consists of all groups acting on trees, with isotropy in $\mathcal H_{i-1}$.
\end{itemize}

The main theorem which will allow us to reduce the Baum-Connes conjecture for any finitely generated subgroup of $SL(2,\C)$ to the groups of integral characteristic will then be:

\begin{theorem}\label{SL2CinH}
Every finitely generated subgroup of $SL(2,\mathbb C)$ lies in $\cup_{i=0}^{\infty}\mathcal H_i$.
\end{theorem}

Indeed, repeated applications of Theorem~\ref{OyonoOyonoTheorem} allow one to reduce the Baum-Connes conjecture for any group in $\mathcal H_i$ to groups in $\mathcal H_0$, that is groups of integral characteristic.

The rest of this section is devoted to the proof of Theorem~\ref{SL2CinH}.

We start with a finitely generated subgroup $\Gamma$ of $SL(2,\C)$ and notice that we can treat it as a subgroup of $SL(2,A)$, where $A$ is a finitely generated subring of $\C$.

The ultimate goal is to show that $\Gamma$ lies somewhere in the sequence $\{\mathcal H_i\}$, defined above.

Now we construct a special hierarchy of subgroups of $\Gamma$, resembling the technics of Alperin and Shalen in~\cite{AlperinShalen} (see Chapter~\ref{chapterMainTheorem} for a detailed discussion; construction here is  closer to the original one in~\cite{AlperinShalen} and treats every term of the diagonal action constructed in~Chapter~\ref{chapterMainTheorem} one-by-one) and show that this hierarchy fits inside the general one, described above. Starting from $\Gamma\subseteq SL(2,A)$, construct a sequence $\{\tilde{\mathcal H}_0,\dots,\tilde{\mathcal H}_{m+1}\}$ of families of subgroups of $\Gamma$ in the following fashion:

\begin{itemize}
\item $\tilde{\mathcal H}_{m+1}=\{\Gamma\}$
\item $\tilde{\mathcal H}_i=\{$ subgroups of $\Gamma$ with the coefficients of the characteristic polynomial of every element in $\cap_{j=i}^m\mathcal O_{\nu_j}\}$ for $i=m,m-1,\dots,0$ 
\end{itemize}
(As before, we use the following notation: $\nu_0,\dots,\nu_m$ are the discrete valuations furnished by Lemma~\ref{lemmaMainThmValuations} for the field of fractions of $A$; $\mathcal O_{\nu_i}$ for $i=0,\dots,m$ are the corresponding rings of integers, and we consider an action of $\Gamma$ and its subgroups on the trees of the equivalence classes of $\mathcal O_{\nu_i}$-lattices.)

It is easy to see that for each $i>0$ any group $G$ in $\tilde{\mathcal H}_i$ acts on the  $\mathcal O_{\nu_{i-1}}$-tree with isotropy in subgroups of conjugates of $SL(2,\mathcal O_{\nu_{i-1}})$, that is in $\tilde{\mathcal H}_{i-1}$. On the other hand, Lemma~\ref{lemmaMainThmValuations} shows that the intersection of all rings of integers $O_{\nu_i}$ is a subring in the algebraic integers, that is to say any group in $\tilde{\mathcal H}_0$ has integral characteristic. Thus we have shown that $\tilde{\mathcal H}_i\subseteq\mathcal H_i$ for any $i=0,1,\dots,m$, whence Theorem~\ref{SL2CinH} for our particular group $\Gamma$.

\subsection{Zariski-dense subgroups}\label{ZariskiDense}

Now we need to prove the Baum-Connes conjecture for subgroups of integral characteristic. In this section we concentrate on integral characteristic subgroups $\Gamma$ of $SL(2,\C)$ whose Zariski closure $G$ is the entire $SL(2,\C)$ (Zariski density here corresponds to the irreducibility of the action on $K^n$ in Section~\ref{sectionIntCharIrrAction}).

The following result is a modification of Lemma~\ref{lemmaIntCharEmbeddingAlpha} and it goes back to Zimmer (cf.~\cite[Lemma 6.1.7]{ZimmerErgodicTheory}).

\begin{lemma}\label{lemmaSL2ZimmerEmbedding}
There exists a faithful representation
\begin{equation*}
\alpha: SL(2,\C)\to GL(4,\C),
\end{equation*}
such that $\alpha(\Gamma)\subset GL(4,\A)$.
\end{lemma}
\begin{proof}
We shall use the construction, similar to the one in the proof of Lemma~\ref{lemmaIntCharEmbeddingAlpha}. Let $f_g$ be a complex-valued map $G\to\C$ defined by
\begin{equation*}
f_g:h\mapsto\tr(gh), \qquad h\in G.
\end{equation*}
Note that $f_{g_1}(h)+f_{g_2}(h)=\tr((g_1+g_2)h)$ and $\lambda f_g(h)=\tr((\lambda g)h)$ for any $g_1, g_2, h\in G$ and $\lambda\in\C$. This allows us to consider $f_g(h)$ as a short-hand notation for $\tr(gh)$ for any $g\in\C G$, $h\in G$. Let
\begin{equation*}
V=\Span_{\C}\left\lbrace f_g \right\rbrace_{g\in\C G}.
\end{equation*}
Since the conditions defining $f_g$ are linear with respect to the entries of $g$, the linear space $V$ has finite dimension. More precisely, its dimension is $4$.

Consider the following action of $G$ on $V$:
\begin{equation}\label{equationSL2ActionZDense}
g.f_h=f_{gh}, \qquad g\in G, h\in\C G.
\end{equation}
Since this action is linear, we have a representation of $G$.

Let
\begin{equation*}
W=\Span_{\C}\left\lbrace f_g \right\rbrace_{g\in\C\Gamma}.
\end{equation*}
This subspace is $\Gamma$-invariant and, since $\Gamma$ is Zariski-dense in $G$, is also $G$-invariant. Thus $W=V$.

Let $g_1, g_2, g_3, g_4\in\Gamma$ be such that $\left\lbrace f_{g_1}, f_{g_2}, f_{g_3}, f_{g_4}\right\rbrace $ is a basis of $V$ (we can arrange this because $V$ is generated by $f_g$ for $g\in\Gamma$). With respect to this basis the action~\eqref{equationSL2ActionZDense} is given by a matrix $(\alpha_{ij}^g)$, such that
\begin{equation}\label{equationSL2ActionMatrix}
g.f_{g_i}(h)=\sum_{j=1}^4\alpha_{ij}^gf_{g_j}(h), \qquad g, h\in G, i=1,\dots,4.
\end{equation}
Thus we obtain a representation $\alpha: G\to GL(4,\C)$.

Similarly to the proof of Lemma~\ref{lemmaIntCharEmbeddingAlpha}, we then show that $\alpha$ is faithful.

If $g\in\Gamma$, \eqref{equationSL2ActionMatrix} means that in particular
\begin{equation*}
\tr(gg_ig_k)=g.f_{g_i}(g_k)=\sum_{j=1}^4\alpha_{ij}^gf_{g_j}(g_k)=\sum_{j=1}^4\alpha_{ij}^g\tr(g_jg_k), \qquad i, k=1,\dots,4.
\end{equation*}
Then $(\alpha_{ij}^g)$ are the solutions of this system of linear equations with algebraic coefficients, and therefore the matrix entries $(\alpha_{ij}^g)$ of the representation $\alpha$ are algebraic.
\end{proof}

Since the conclusion of Lemma~\ref{lemmaSL2ZimmerEmbedding} is the same as the one of Lemma~\ref{lemmaIntCharEmbeddingAlpha}, for any other Zariski--dense subgroup $\tilde\Gamma$ of integral characteristic its image $\alpha(\tilde\Gamma)$ is conjugate in $GL(4,\C)$ to a subgroup whose matrix entries belong to a finitely generated subring of $\A$. In the discussion which follows we continue to work with $\Gamma$.

By the Primitive Element Theorem $\alpha(\Gamma)\subset GL(4,K)$, with $[K:\Q]<\infty$. Take $H=\alpha(G)$. We see that $\alpha(\Gamma)\subseteq H\cap GL(4,K)$,
and both $\alpha(\Gamma)$ and $H\cap GL(4,K)$ are Zariski--dense in $H$, thus they are both defined over $K$ by~\cite[Proposition 3.1.8]{ZimmerErgodicTheory}. This means $\alpha(\Gamma)$ is locally isomorphic to $SL(2,K)$ (see~\cite[Theorem~7]{Zimmer}), and we can apply Theorem~\ref{theoremIntCharIrr} to $\alpha(\Gamma)$ and obtain a proper action. Via a finite-dimensional version of the result of Higson and Kasparov in~\cite{HigsonKasparov}, the Baum-Connes conjecture for our group follows. 

\subsection{Zariski--non--dense subgroups}\label{ZariskiNonDense}

Finally, we discuss the case of Zariski--non--dense subgroups of integral characteristic.

Let us start with some preliminary remarks on algebraic Lie groups.
Suppose $G$ is a Zariski--closed proper subgroup of $SL(2,\C)$. We write $G_0$ for the Zariski--connected component of the unit of $G$. It is well known that $G_0$ is a normal subgroup of $G$ of finite index \cite[I.1.2]{Borel}.

That is to say, $G$ is almost connected in this case, and the Baum-Connes conjecture for almost connected groups is proven in~\cite{ChabertEchterhoffNest} and~\cite{Lafforgue}. This general result requires quite heavy $KK$-theory technique and is proven for the Baum-Connes conjecture with trivial coefficients only. Here we shall present a more elementary argument for the version of the conjecture with coefficients, by explicitly describing the structure of the groups we are dealing with.

Since for algebraic groups the notions of connected and irreducible components coincide~\cite[AG.17.2]{Borel}, $G_0$ is abelian if and only if its Lie algebra is commutative \cite[IV.4.3]{Hochschild}. 
Since $G$ is a proper subgroup of $SL(2,\C)$, its dimension is strictly less than $3$. In the subsections below we shall address each dimension case separately.

\subsection*{Dimension $0$}

In this case $\dim G_0=0$ as well, and, since $G_0$ is connected, we conclude that it is trivial. The group $G$ itself, being a finite extension of $G_0$, is finite, whence the Baum-Connes conjecture for $G$ holds trivially.

\subsection*{Dimension $1$}

The Lie algebra of $G$ (and $G_0$ as well) has to be $1$-dimensional. In particular, it has to be commutative, thus $G_0$ is abelian. There are only two (up to conjugacy) connected abelian $1$-dimensional groups, namely $\left\{\left. \m{a}{0}{0}{a^{-1}}\right| a\in\C^{\times}\right\}$ and $\left\{\left. \m{1}{b}{0}{1}\right| b\in\C\right\}$. We shall treat them separately.

Suppose $G_0=\left\{\m{a}{0}{0}{a^{-1}}\right\}$ (up to conjugacy). Then, since $G_0$ is normal in $G$, any conjugate of $\m{a}{0}{0}{a^{-1}}$ by any element in $G$, say $\m{g_{11}}{g_{12}}{g_{21}}{g_{22}}$, has to have the same diagonal form:
\begin{multline*}
\m{g_{11}}{g_{12}}{g_{21}}{g_{22}}\m{a}{0}{0}{a^{-1}}\m{g_{11}}{g_{12}}{g_{21}}{g_{22}}^{-1}=
\m{*
}{(a^{-1}-a)g_{11}g_{12}}{(a-a^{-1})g_{21}g_{22}}{*
}=
\m{b}{0}{0}{b^{-1}}.
\end{multline*}
This means that $g_{11}g_{12}=0$ and $g_{21}g_{22}=0$. To satisfy the first condition, we need to take either $g_{11}$ to be zero or $g_{12}$ to be zero. Thus $G$ may contain only matrices with zeros on the diagonal, or off the diagonal:
\begin{equation*}
G\subseteq\left\{\m{a}{0}{0}{a^{-1}}, \m{0}{a}{-a^{-1}}{0}\right\}=H.
\end{equation*}
This group $H$ is amenable, and, modifying Theorem~\ref{theoremOyonoOyono}, it is possible to show (see~\cite[Theorems 5.18 and 5.23]{MislinValette}) that any countable subgroup of $H$ satisfies the Baum-Connes conjecture with coefficients\footnote{Theorem 5.18 in~\cite{MislinValette} shows that $H$ satisfies the Baum-Connes conjecture with trivial coefficients, while Theorem 5.23 proves that any countable subgroup of such group satisfies the conjecture with arbitrary coefficients.}. 

Now suppose $G_0=\left\{\m{1}{b}{0}{1}\right\}$ (again, up to conjugacy). Since $G_0$ is normal in $G$, any conjugate of $\m{1}{b}{0}{1}$ by an arbitrary element $\m{g_{11}}{g_{12}}{g_{21}}{g_{22}}$ element of $G$ should have the same form:
\begin{equation*}
\m{g_{11}}{g_{12}}{g_{21}}{g_{22}}\m{1}{b}{0}{1}\m{g_{11}}{g_{12}}{g_{21}}{g_{22}}^{-1}=\m{*}{*}{-bg_{21}^2}{*}=\m{*}{*}{0}{*}.
\end{equation*}
Thus we have $g_{21}=0$, which means that $G$ contains only matrices of the form $\m{a}{b}{0}{a^{-1}}$, and we shall discuss this group in the following subsection.

\subsection*{Dimension $2$}

Let $H$ denote the subgroup of $SL(2,\mathbb C)$ consisting of all the matrices of the form $\m{a}{b}{0}{a^{-1}}$, where $a\in\C^{\times}, b\in\C$. Note that its Lie algebra consists of the matrices of the form $\m{a}{b}{0}{-a}$, $a, b\in\C$.
\begin{lemma}\label{MaximalityOfH}
Any subgroup $K$ of $SL(2,\mathbb C)$, which includes $H$ and some element not in $H$, coincides with the whole group  $SL(2,\mathbb C)$.
\end{lemma}
\begin{proof}
Suppose $K$ contains some element \m{g_{11}}{g_{12}}{g_{21}}{g_{22}} with $g_{21}\ne0$. Then we can multiply this element by an element \m{g_{21}^{-1}}{-g_{22}}{0}{g_{21}} of $H$ on the right to get \m{g_{11}g_{21}^{-1}}{-1}{1}{0}.

Now for any complex numbers $a$ and $b$ with $a\ne0$ we can multiply \m{g_{11}g_{21}^{-1}}{-1}{1}{0} on the left by \m{a}{b-ag_{11}g_{21}^{-1}}{0}{a^{-1}} to get \m{b}{-a}{a^{-1}}{0}. Since $K$ is a group, it ought to contain all inverses as well, in particular \m{0}{a}{-a^{-1}}{b}. This means that $K$ contains all matrices of determinant $1$ with $0$ in the upper left corner.

Finally, let us take an arbitrary element of $SL(2,\mathbb C)$, say \m{s_{11}}{s_{12}}{s_{21}}{s_{22}}. Since we already know that all matrices with  $s_{21}=0$ belong to $H$, and therefore to $K$, the essential part of the argument is to show that any such matrix with $s_{21}\ne0$ belongs to $K$. The identity
\begin{equation*}
\m{s_{11}}{s_{12}}{s_{21}}{s_{22}}=\m{1}{s_{11}s_{21}^{-1}}{0}{1}\m{0}{-s_{21}^{-1}}{s_{21}}{s_{22}}
\end{equation*}
completes the proof, since all matrices on the right-hand-side belong to $K$.
\end{proof}

Now we provide some technical results about Lie subalgebras of $\mathfrak{sl}(2,\C)$.
\begin{lemma}
For any 2-dimensional noncommutative Lie algebra there exists a basis $\{x_1, x_2\}$ with multiplication table $[x_1x_2]=x_1$.
\end{lemma}
\begin{proof}
Consider some basis $\{e_1, e_2\}$ of such an algebra, and suppose $[e_1e_2]=a_1e_1+a_2e_2$. Since the algebra is noncommutative, at least one of the coefficients $a_1$ and $a_2$ is not zero, let us say $a_1\ne0$, for definiteness. It is easy to see that two elements $x_1=a_1e_1+a_2e_2$ and $x_2=a_1^{-1}e_2$ are linearly independent, so that they constitute a basis as well. With respect to such basis the multiplication is given by $[x_1x_2]=x_1$.
\end{proof}
\begin{lemma}\label{OnlyOne}
The Lie algebra $\mathfrak{sl}(2,\C)$ contains only one (up to conjugation) $2$-dimensional Lie subalgebra, namely $\left\{\left.\m{a}{b}{0}{-a}\right| a, b\in\C\right\}$.
\end{lemma}
\begin{proof}
Suppose we have a $2$-dimensional noncommutative Lie subalgebra $\mathfrak h$ of $\mathfrak{sl}(2,\C)$. Let $\{x_1, x_2\}$ denote the basis of $\mathfrak h$ constructed in the lemma above. A priori there could be two possibilities: both eigenvalues of the matrix $x_2$ coincide (and therefore are zeros) or they are distinct. In the first case $x_2$ is conjugate to its Jordan form, namely $\m{0}{1}{0}{0}$, and if $x_1$ after same conjugation has the form $\m{a}{b}{c}{d}$, then the multiplication condition $[x_1x_2]=x_1$ is
\begin{equation*}
\left[\m{a}{b}{c}{d},\m{0}{1}{0}{0}\right]=\m{-c}{a-b}{0}{c}=\m{a}{b}{c}{d},
\end{equation*}
from where we conclude that $a=b=c=d=0$, which means $x_1=0$ and therefore can not serve as basis element, so that the case where both eigenvalues of the matrix $x_2$ coincide can not happen. Now suppose that the eigenvalues of $x_2$ are distinct, say $\lambda$ and $-\lambda$. Conjugating $x_1$ and $x_2$, we write the multiplication condition as
\begin{equation*}
\left[\m{a}{b}{c}{d},\m{\lambda}{0}{0}{-\lambda}\right]=\m{0}{-2b\lambda}{2c\lambda}{0}=\m{a}{b}{c}{d},
\end{equation*}
so that we have $a=d=0$ and $2c\lambda=c$, $-2b\lambda=b$. We are looking for solutions with at least one of the coefficients $b$ and $c$ being non-zero, therefore we end up with two possibilities:
\begin{enumerate}
\item $b=0\ne c$, $\lambda=\frac12$
\item $c=0\ne b$, $\lambda=-\frac12$
\end{enumerate}
Thus any non-commutative $2$-dimensional Lie subalgebra of $\mathfrak{sl}(2,\mathbb C)$ is conjugate-equivalent to $\mathfrak h_1=\mathbb C\m{0}{0}{1}{0}\oplus\mathbb C\m{\frac12}{0}{0}{-\frac12}$ or $\mathfrak h_2=\mathbb C\m{0}{1}{0}{0}\oplus\mathbb C\m{-\frac12}{0}{0}{\frac12}$. Finally, $\mathfrak h_1$ and $\mathfrak h_2$ are conjugate to each other via the matrix \m{0}{1}{1}{0}. By scaling the second matrix in $\mathfrak h_2$ , we obtain the representation $\left\{\left. \m{a}{b}{0}{-a} \right| a, b\in\mathbb C\right\}$.

Now we show that $\mathfrak{sl}(2,\C)$ does not contain any commutative $2$-dimensional subalgebras. Suppose one such exists and has a basis $\{x, y\}$. Conjugating by some matrix, we can put $y$ into Jordan form, and let $x$ be represented by \m{x_{11}}{x_{12}}{x_{21}}{-x_{11}} under the same conjugation. We have two possibilities: the eigenvalues of the matrix, representing $y$, coincide (and therefore are zeros) or they are distinct, and by scaling the matrix we assume that they are $1$ and $-1$. In the first case the commutativity condition can be written as
\begin{equation*}
\m{x_{11}}{x_{12}}{x_{21}}{-x_{11}}\m{0}{1}{0}{0}=\m{0}{1}{0}{0}\m{x_{11}}{x_{12}}{x_{21}}{-x_{11}},
\end{equation*}
which leads to $x_{11}=x_{21}=0$, so that $x$ is a scalar multiple of $y$, and this can not happen. In the second case we have
\begin{equation*}
\m{x_{11}}{x_{12}}{x_{21}}{-x_{11}}\m{1}{0}{0}{-1}=\m{1}{0}{0}{-1}\m{x_{11}}{x_{12}}{x_{21}}{-x_{11}},
\end{equation*}
this means $x_{12}=x_{21}=0$, and again we have a contradiction with linear independence of $x$ and $y$.
\end{proof}

Getting back to the group $G_0$, we see that Lemma~\ref{OnlyOne} describes the Lie algebra of $G_0$, up to conjugacy. Therefore $G_0$ and $H$ are conjugate to each other. 
Lemma~\ref{MaximalityOfH} suggests that there are no proper subgroups of $SL(2,\C)$, larger than $H$, therefore we conclude that $G=G_0$.

Finally, $G$ is a semidirect product
\begin{equation*}
\left\{\left. \begin{pmatrix}1&b\\0&1\end{pmatrix} \right| b\in\C\right\}\rtimes\left\{\left. \begin{pmatrix}a&0\\0&a^{-1}\end{pmatrix} \right| a\in\C^{\times}\right\}
\end{equation*}
of two abelian groups, hence it is amenable, and we conclude that any finitely generated subgroup of $G$ satisfies the Baum-Connes conjecture with coefficients by applying~\cite[Theorem 5.23]{MislinValette}.

Now Theorem~\ref{theoremMain2} has been proven completely.